\documentclass[12pt]{article}

\usepackage{
amssymb,amsmath}

\def\l{\lambda}
\def\xx{{\cal X}}
\def\ww{{\cal W}}
\def\tG{{\widetilde G}}
\def\rmd{{\rm d}}
\def\var{\text{\rm Var\,}}
\def\r{\rho}
\def\s{\sigma}
\def\a{\alpha}
\def\b{\beta}
\def\g{\gamma}

\def\d{\delta}

\def\e{\varepsilon}
\def\h{\eta}
\def\th{\theta}
\def\m{\mu}

\def\t{\tau}
\def\f{\varphi}

\def\ps{\psi}
\def\giv{\,|\,}
\def\half{{\tfrac12}}
\def\shalf{{\scriptstyle {1 \over 2}}}
\def\quarter{{\textstyle {1 \over 4}}}
\def\nin{\noindent}

\def\msk{\medskip}
\def\mm{{\tt m}}
\def\tinu{{\tilde\nu}}
\def\tS{{\widetilde S}}
\def\tPhi{{\widetilde \Phi}}
\def\tL{{\widetilde L}}
\def\tPsi{{\widetilde \Psi}}
\def\ignore#1{}

\def\k{\varkappa}
\def\Le{\ \le\ }
\def\Ge{\ \ge\ }
\def\Eq{\ =\ }
\def\KK{{\cal K}}
\def\dt{\,{\rm d}t}
\def\ds{\,{\rm d}s}
\def\dx{\,{\rm d}x}
\def\dy{\,{\rm d}y}
\def\du{\,{\rm d}u}
\def\dv{\,{\rm d}v}
\def\dw{\,{\rm d}w}
\def\dz{\,{\rm d}z}
\def\kvi{(\k\vee1)}
\def\Refrm#1{{\rm (\ref{#1})}}

\def\hB{{\widehat B}}
\def\ut{^{(2)}}
\def\hX{{\widehat X}}
\def\kk{{\cal K}}
\def\nn{{\cal N}}
\def\ui{^{(1)}}
\def\ut{^{(2)}}
\def\uh{^{(3)}}
\def\tY{{\widetilde Y}}
\def\thrhalf{\textstyle{3\over 2}}
\def\shalf{\scriptstyle{1\over 2}}

\newcommand{\eqs}{\begin{eqnarray*}}
\newcommand{\ens}{\end{eqnarray*}}

\newcommand{\eqa}{\begin{eqnarray}}
\newcommand{\ena}{\end{eqnarray}}
\newcommand{\eq}{\begin{equation}}
\newcommand{\en}{\end{equation}}

\def\numberlikeadb{\global\def\theequation{\thesection.\arabic{equation}}}
\numberlikeadb

\newtheorem{theorem}{Theorem}[section]
\newtheorem{lemma}[theorem]{Lemma}
\newtheorem{corollary}[theorem]{Corollary}
\newtheorem{proposition}[theorem]{Proposition}

\def\pr{{\mathbb P}}
\def\ex{{\mathbb E}}
\def\re{{\mathbb R}}

\def\sqln{\sqrt{\log n}}
\def\sqLn{\sqrt{L(n)}}
\def\sqpsn{\sqrt{\ps_n}}
\def\tnT{{\t_n\wedge T}}
\def\psnT{{\ps_n\wedge T}}

\def\hZ{{\widehat Z}}
 \newcommand{\tnu}{{\nu}}

\newcommand{\hPhi}{{\Phi}}

\newcommand{\hPsi}{{\Psi}}

\def\endpf{\hfill $\Box$ }

\def\Kni{{\cal K}_n(\infty)}
\def\ul{^{(l)}}
\def\thalf{\tfrac12}

\def\ff{{\cal F}}
\def\Ref#1{(\ref{#1})}
\def\Le{\ \le\ }

\def\Blm{\left|}
\def\Brm{\right|}
\def\Blb{\left\{}
\def\Brb{\right\}}
\def\Bl{\left(}
\def\Br{\right)}
\def\nti{n\to\infty}
\def\lnti{\lim_{\nti}}
\def\law{{\cal L}}

\def\bone{{\bf 1}}
\def\non{\nonumber}

\begin{document}

\title{
Regenerative Compositions in the Case of Slow Variation}
\author{A.D.~Barbour\thanks{Angewandte Mathematik, Winterthurerstrasse~190,
CH--8057 Z\"urich, Switzerland; supported in part by Schweizer Nationalfonds
Projekt Nr.~20-107935/1} ~and A.V.~Gnedin\thanks{
Mathematisch Instituut, PO Box 80010, 3508 TA Utrecht, 
Nederland}\\
Universit\"at Z\"urich and Universiteit Utrecht}
\date{}
\maketitle

\begin{abstract}
\noindent
For $S$ a subordinator and $\Pi_n$ an independent Poisson process of intensity $ne^{-x}, x>0,$ 
we are interested in the number $K_n$ of gaps in the range  of $S$ that are  hit by at 
least one point of $\Pi_n$. Extending previous studies in \cite{Bernoulli, GPYI, GPYII}
we focus on the  case when the tail of the L{\'e}vy measure of $S$ is slowly varying.
We view $K_n$ as the terminal value of a random process ${\cal K}_n$, and 
provide an asymptotic analysis of the fluctuations of ${\cal K}_n$, as $n\to\infty$, for a
wide spectrum of situations.  

\end{abstract}

\setcounter{section}0

\section{Introduction}

Let $S = (S_t, t\geq 0)$ be an increasing L{\'e}vy process (subordinator) 
with $S_0=0$, zero drift and no killing.
The closed range ${\cal R}$ of~$S$ has zero Lebesgue measure, and 
defines a random division of the complement set ${\mathbb R}_+\setminus {\cal R}$ 
into open interval components, referred to as {\it gaps\/}.  
In this paper, we are concerned with the 
distribution of the number $K_n$ of gaps hit 
by at least one point of an independent Poisson process~$\Pi_n$ with the 
inhomogeneous rate $n e^{-x}$, $x>0$, where~$n$ is a large parameter. 
We actually go into more detail. We view~$K_n$ as the terminal value~$\Kni$ 
of the increasing process ${\cal K}_n=({\cal K}_n(T),\, T\geq 0)$,
where ${\cal K}_n(T)$ is defined to be the number of jumps of the 
subordinator within $[0,T]$ which cover one or more Poisson points; 
that is, ${\cal K}_n(T)$ counts all instants 
$t\in [0,T]$ which satisfy $\Pi_n\cap \,]S_{t-}, S_t[\,\neq \emptyset$.
Our aim is to describe the random fluctuations of the process~$\kk_n$.

\par The general motivation for the setting stems from the study of the 
number of blocks in a random decomposable 
combinatorial structure, which in the case under focus  is a 
{\it composition\/} (ordered partition) ${\cal C}_n$ of some integer.
Viewing ${\cal C}_n$ as distribution of some number of balls 
in some collection of boxes, each gap may be interpreted as a box, which is
hit by a particular ball with probability equal to the exponential measure of 
the gap. The parameter~$n$ controls the total number of balls, which 
is a Poisson variable, and the  composition ${\cal C}_n$ 
of this random number is defined  as the consecutive record of nonzero occupancy 
numbers, in the natural ordering of the gaps.
Our study fits into the recent theory of sampling models 
called {\it regenerative composition structures\/}, which have  
a distinguished Markovian property resulting from the renewal
features of  ${\cal R}$ combined with that
of the exponential distribution \cite{RCS, RPS}.
Concretely, the regeneration property of ${\cal C}_n$ means that, for each 
$t>0$, 
conditionally
given the value $s=S_t$, 
the partial compositions 
appearing within $[0, s]$ and $[s,\infty[$ are independent and the latter has the same 
distribution as the composition ${\cal C}_{ne^{-s}}$.

\par The distribution of~$S$ is completely determined by a L{\'e}vy measure $\nu_0$ on~$\re_+$, which
describes the intensity of the jumps of different sizes, and the
behaviour of~$\kk_n$ depends very much on the form of~$\nu_0$.  
Qualitatively different modes of behaviour are known \cite{Bernoulli, GPYI, GPYII}.

\par For~$\nu_0$  a finite measure, $S$ is a compound Poisson process.
Under mild additional assumptions, 
the two central moments are of the order of $\log n$ and $K_n$ is asymptotically normal.
In this situation, the methods of renewal theory are adequate, 
since
the process ${\cal K}_n(T)$ essentially coincides with the process of jump epochs of $S$ for $T<\log n$, while the 
contribution of larger times $T>\log n$ to $K_n$ is negligible,
 see \cite{Bernoulli}.
In particular, when
$\nu_0$ is an exponential distribution,
the induced composition follows 
the poissonised (ordered) Ewens sampling formula, in which case  much finer results on $K_n$
are available by combinatorial methods
\cite{ABT, CSP}.

\par If~$\nu_0$ is infinite, and its tail $N_0(x):=\nu_0[x,\infty[\, $ is such that
$N(1/y)$ is regularly varying as $y\to\infty$ with exponent~$\a$ 
(here and henceforth this means regular variation with $0<\alpha\leq 1$), then 
${\mathbb E}K_n$ is also regularly varying with the same exponent and $K_n/{\mathbb E}K_n$ approaches a nondegenerate 
limit, which is not gaussian. The moments of
${\cal K}_n(T)$ are then of the same order of magnitude as that of $K_n$, for each fixed $T$,
see \cite{GPYI}.

\par Between these two possibilities
lies the setting in which $N_0(1/y)$ is slowly varying as $y\to\infty$, but the L{\'e}vy measure is infinite,
{\it i.e.}
$\lim_{y\to\infty} N_0(1/y) = \infty$.
 Here, the special case with 
$N_0(1/y)\sim c \log y$ has been studied in some detail. For these
 {\it gamma-like\/} subordinators, the proper formats for the two central moments of $K_n$ are 
$\log^2 n$ and $\log^3 n$, respectively, and the limiting distribution is again normal, see \cite{GPYII}.

\par In this paper, we treat the case of slowly varying~$N_0$ in greater
generality.  As might be expected of a transitional r\'egime between
the finite and the regularly varying cases, there is a further wealth of
possible modes of behaviour, and the discussion reveals how these 
are related to the time scales over which the significant variation
in~$\kk_n$ occurs. Our argument leading to a functional central limit theorem 
is very different from that in~
\cite{Bernoulli, GPYII},
and is based on the observation that, to first order, the fluctuations
of the counting process~$\kk_n$ are dominated by those of its
compensator~$A_n$, defined in Proposition~\ref{A-n-def}.
The explicit representation
of the random process~$A_n$
 makes it possible to find approximations by rather direct
arguments, and under relatively mild conditions.
These are broadly speaking of two kinds.  The first is expressed
in Assumption~A2, which puts
a mild restriction on the way in which a certain transform 
$L$ 
of the
measure~$\nu_0$ can vary locally as a function of its parameter. 
Conditions of the second kind, appearing in different forms in 
\Ref{L-cond1}, \Ref{L-cond2},
\Ref{L-cond4} and~\Ref{L-cond6}, limit
the global variability of $L$.

\par Our analysis  of subordinators with slowly varying $N_0$
distinguishes three basic modes. 
In the case of {\it moderate growth}, which includes 
the subordinators with logarithmic asymptotics
$N_0(1/y)\asymp (\log y)^\beta$, $\beta>0$, including the gamma-like 
subordinators studied in \cite{GPYII}, the random 
fluctuations of ${\cal K}_n(T)$ occur more or less evenly on the 
scale $T=v\log n$ in $v\in [0,1]$. In the case of {\it fast growth\/}, 
well exemplified by $N_0(1/y)\asymp \exp(\log^\beta y)$, $0<\beta<1$,
almost everything happens at times of order~$L(n)$, and~$L(n)$ is of 
smaller order than $\log n$. The third case is that of {\it slow growth\/}, as 
for example $N_0(1/y)\asymp \log \log y$, when significant contributions
to the random fluctuations of~$K_n$ are only made at times very close 
to $\log n$, just as in the compound Poisson case \cite{Bernoulli}.

\vskip0.5cm
\noindent
{\bf Notation.} We  use $\lambda$, $\lambda_j$ for positive 
constants whose value is not important and may depend on the  context.
The asymptotic relation $a_n\asymp b_n$ means that $a_n=O(b_n)$ and $b_n=O(a_n)$, 
while $a_n\gg b_n$ means that $b_n=o(a_n)$. Asymptotic relations like $X_n\sim Y_n$ or $X_n\asymp Y_n$ 
for random quantities mean that they hold with probability one, 
unless otherwise specified.

\section{The basic setting}
\subsection{Laplace exponents and the compensator}

The L{\'e}vy measure~$\nu_0$ is uniqely determined by the
Laplace exponent~$\Phi_0$, defined for $m\geq 0$ by
$$
  \Phi_0(m) \Eq \int_0^\infty (1-e^{-mx})\,\nu_0({\rm d}x)
	  \Eq m \int_0^\infty e^{-mx}N_0(x)\,{\rm d}x\,;
$$
note that~$\nu_0$ must satisfy $\Phi_0(1)<\infty$. The distribution of 
the subordinator
is determined by the L{\'e}vy-Khintchine formula for the Laplace transform 
\eq\label{LK}
\ex\,e^{-nS_t}=e^{-t\Phi_0(n)}\,,~~~~n\geq 0,
\en
see \cite{Bertoin} as a general reference on the L{\'e}vy processes and see
\cite{BertoinSub} especially for subordinators.

\par The function $\Phi_0$ can be extended to an analytic function in the 
right half-plane, and hence, by M{\"u}ntz's theorem, 
$\Phi_0$ can  be uniquely extrapolated from the values 
$\Phi_0(m)$, $m=1,2,\ldots$; these also determine the 
poissonised version of~$\Phi_0$, defined either by the series
\eq\label{hPhi}
  \hPhi(n)\ :=\ e^{-n}\sum_{m=1}^\infty \frac{n^m}{m!}\Phi_0(m)\,,
\en
or by the integral
\eq\label{hPhiIn}
  \hPhi(n)\Eq\int_0^\infty (1-e^{-n(1-e^{-x})})\,\nu_0({\rm d}x).
\en
This latter transform is particularly useful to us, since it
appears naturally in the definition of the compensator~$A_n$
of the counting process~$\kk_n$.

\begin{proposition}\label{A-n-def}
With respect to the filtration
$(\ff_{T,n},\,T\ge0)$, defined by
\[
  \ff_{T,n}\ :=\ \s\Blb S_t,\,0\le t\le T;\ \Pi_n|_{[0,S_T]}\Brb,
\]
the compensator of ${\cal K}_n$ 
is the increasing process~$A_n$ given by the formula
\eq\label{compensator-def}
  A_n(T)\ :=\ \int_0^T \,\hPhi(n e^{-S_t})\,{\rm d}t\,,\qquad T\in [0,\infty]\,.
\en
\end{proposition}

\noindent
\begin{proof} The subordinator gains an increment within $[x, x+{\rm d}x]$ 
at rate $\nu_0({\rm d}x)$. On the other hand, $\Pi_n$ hits $[S_{t-}, S_{t-}+x]$ 
with probability $1-\exp({-n e^{-s}(1-e^{-x})})$ (where $s=S_{t-}$),
because the  number of atoms in $[s\,,s+x]$ has Poisson distribution 
with mean
$$
  n \int_s^{s+x} e^{-u}\,{\rm d}u= n e^{-s}(1-e^{-x})\,.
$$
Integrating over $x$ yields the derivative 
${\rm d}A_n(t)/{\rm d}t=\hPhi(ne^{-S_t})$.
\endpf\end{proof}

\bigskip
\ignore{
\noindent Since the definition of~$A_n$ does not involve the process~$Y_n$,
it follows that~${\cal K}_n$ is actually a Cox process, with directing measure
determined by the subordinator~$S$.  In particular, we have
\eq\label{ADB:mean-var}
  \ex K_n \Eq \ex A_n(\infty)\quad {\rm and}\quad \var K_n 
     \Eq \var A_n(\infty) + \ex A_n(\infty)
\en
(c.f. Daley and Vere--Jones~\cite{DVJ}, Theorem 13.4.I and p.~263).
}     

\par
We further assume that 
$${\mathbb E}S_t=t ~~~~{\rm and}~~~~ \var S_t = t\s^2\,.$$
The former is the same as 
\eq\label{drift-one}
  {\mathbb E}S_1=\Phi_0'(0)=\int_0^\infty x\,\nu_0({\rm d}x)
	  = \int_0^\infty N_0(x)\,{\rm d}x=1\,.
\en
This can always be achieved by a linear time-scaling, which does not affect 
the range ${\cal R}$.  
A consequence of the assumption is that $N_0(x)$ is a probability density.
It then follows for all $m\ge0$ that
\eq\label{Phi-bnd-1}
  \Phi_0(m) \Le m\int_0^\infty N_0(x)\,{\rm d}x \Eq m\quad\mbox{and}\quad
	\hPhi(m) \Le m,
\en
this last from~\Ref{hPhi}.	

\par With this scaling,
 $N_0(x)$ is the density of a {\it delay} variable.
If $X$ has this density and is independent of~$S$, then
the process $(X+S_t,\, t\geq 0)$ is a {\it stationary} subordinator, 
in the sense that 
its closed range $X+{\cal R}$ may be extended to a random subset 
of ${\mathbb R}$ invariant under all translations. 
In particular, $X+S_t$ 
has the same overshoot distribution 
at every level $s\geq 0$.

\par It will be convenient to define all the Poisson processes 
$(\Pi_n, \,n\geq 0)$ (which are independent of the subordinator~$S$)
consistently on the same probability space. 
To this end, we take an inhomogeneous planar Poisson point 
process on ${\mathbb R}^2_+$ with intensity measure $e^{-y}\dy\,{\rm d}n$,
 and we introduce
$\Pi_n$ as the projection on the $y$-axis of the planar  
process restricted to the strip $[0,\infty]\times[0,n]$.
In this setting, the compositions ${\cal C}_n$ are defined consistently for all 
$n\geq 0$: a decrease in $n$ has the effect of thinning, {\it i.e.} 
removing some balls from the boxes; while as  $n$ increases  more Poisson atoms are added,
hence ${\cal K}_n(T)$ and $K_n$ are nondecreasing in $n$.
Thus, in principle, our setting is 3-dimensional, with three parameters $n,t,s$ meaning the intensity, the time
and the range of subordinator.

\par For our analysis of~$A_n$, it is also convenient
to note that we can truncate the integral~\Ref{compensator-def}, 
which defines~$A_n(T)$, at the first passage time
\eq\label{tau-n-def}
  \tau_n\ :=\ \min\{t\colon S_t\geq \log n\}, 
\en
with little loss.
\begin{lemma}\label{tau-truncation} 
The jumps of $S$ after $\tau_n$ make only a bounded contribution to ${\cal K}_n(T)$, 
uniformly in $T\leq \infty$, and for $\psi>0$
\[
  \pr\left[ \sup_{T\ge0}|A_n(T) - A_n(T\wedge\t_n)| > \ps\right] \Le
	  \ps^{-1}\int_0^\infty \ex\Bl e^{-S_t} \Br\,{\rm d}t.
\]		
\end{lemma}

\begin{proof}$K_n - {\cal K}_n(\tau_n)$ cannot exceed 
the number of atoms of $\Pi_n\cap [\log n, \infty[\,$, which is Poisson 
distributed with mean~$1$.
To estimate the contribution to the compensator, recalling~\Ref{Phi-bnd-1},
we have 
$$
  \int_{\tau_n}^\infty \hPhi(ne^{-S_t})\,{\rm d} t
  \ <\ n \int_{\tau_n}^\infty e^{-S_t}\,{\rm d} t
  \Eq n e^{-S_{\tau_n}}\int_{0}^\infty  e^{-S_u'}\,{\rm d}u
	\Le \int_{0}^\infty e^{-S_u'}\,{\rm d}u\,,
$$
where $S'$ defined by $S'_u := S_{\t_n+u} - S_{\t_n}$, $u\ge0$,
has the same distribution as~$S$. Markov's inequality completes
the proof.
\endpf\end{proof}

\subsection{Slow variation}
Our aim is to investigate the 
process ${\cal K}_n$
in the intermediate setting, between that
in which the tail~$N_0$ of~$\nu_0$ is regularly varying 
(with exponent $0<\alpha\leq 1$), and that
in which it is bounded.  We therefore assume that
\eq\label{slvar}
  \phantom{}\quad{\bf Assumption\  A1}\colon
    N_0(1/y)
{\rm~ is~ an~unbounded~function~ of~ slow~ variation~ for~} y\to\infty.
\en
\noindent
The condition can be equally stated in terms of $\Phi_0$, because
by the Abel--Tauber theorem \cite{Feller}:
$$
  \Phi_0(m)\sim N_0 (1/m)\,,~~~~{\rm as~~} m\to\infty\,.
$$

\par In what follows, we prefer to work in terms of~$\hPhi$, because
of~\eqref{compensator-def}, so that A1 is then more naturally expressed
in the equivalent form: $\hPhi(m)$ is unbounded and slowly varying
at infinity.  
While this equivalence is more or less  clear from~(\ref{hPhi}),  
it is useful for Section~\ref{key} to have a better idea 
of how close the functions $\hPhi$ and~$\Phi_0$ are to each other. 
To this end, we define 
a measure $\tnu$ on $[0,1]$ as the
pushforward of $\nu_0$ under the change of variable $x\to 1-e^{-x}$. 
Then we have
$$
  \Phi_0(m)=\int_0^1 (1-(1-x)^m)\,\tnu({\rm d}x)
	  = m\int_0^1 (1-x)^{m-1}\, N(x) \,{\rm d}x\,,
$$
where $N(x):=\tnu[x,1]$, and
$$
  \hPhi(m)=\int_0^1 (1-e^{-mx})\,\tnu({\rm d}x)
	  = m \int_0^1 e^{-mx}\, N(x) \,{\rm d}x\,.
$$
So the substitution transforms a Laplace integral into a Mellin integral,
while $\hPhi$ assumes the conventional form of a   Laplace exponent
(hence $\hPhi$ also corresponds to some subordinator, whose jump-sizes 
do not exceed $1$). We can now use these representations to prove
the following lemma.

\begin{lemma}\label{2phis} 
Under  Assumption {\rm A1}, we have
$$
  |\Phi_0\ul(m)-\hPhi\ul(m)|=o\left( \frac{\hPhi (m)}{m^{l+1}}\right)\,,
   \quad l\ge 0,
$$
where $f\ul$ denotes the $l$th derivative of~$f$.
\end{lemma}

\noindent{\bf Proof.} Because $N_0(1/y)$ is slowly varying as $y\to\infty$,
the same is true of~$N(1/y)$, and $N(1/m) \sim \hPhi(m)$.  Hence we 
immediately have
\eq\label{Phi-deriv-bnd}
  \hPhi\ul(m) \Eq O\Bl \frac{\hPhi(m)}{m^l} \Br,\quad l\ge 0.
\en
We now define
\[
  D_0(m) \ :=\ m^{-1}\{\hPhi(m) - \Phi_0(m)\}
	  = \int_0^1\{e^{-mx} - e^{(m-1)\log(1-x)}\}\,N(x)\dx.
\]		
Since 
\eqs
  \int_0^1 \{e^{-mx} - e^{m\log(1-x)}\}\,N(x)\dx
	  &=& \int_0^1 e^{-mx}(1 - \exp\{m[\log(1-x)+x]\})\,N(x)\dx\\
		&\sim& \int_0^1 \thalf mx^2 e^{-mx}	N(x)\dx \ \sim\ m^{-2}\hPhi(m),
\ens
and
\[
  \int_0^1 e^{m\log(1-x)}\{1 - (1-x)^{-1}\}\,N(x)\dx
	  \ \sim\ -\int_0^1 xe^{-mx} N(x)\dx \ \sim\ -m^{-2}\hPhi(m),
\]
it then follows that $|D_0(m)| = o(m^{-2}\hPhi(m))$.

For the $l$th derivative $D_l$ of~$D_0$, we similarly have
\[
  D_l(m) = \int_0^1\{(-x)^le^{-mx} 
	  - \{\log(1-x)\}^l e^{(m-1)\log(1-x)}\}\,N(x)\dx.
\]
Once again, since
\eqs
  \lefteqn{\int_0^1 (-x)^l\{e^{-mx} - e^{m\log(1-x)}\}\,N(x)\dx}\\
	  &&\ =\ (-1)^l\int_0^1 x^le^{-mx}(1 - \exp\{m[\log(1-x)+x]\})\,N(x)\dx\\
	&&\ \sim\ (-1)^l\int_0^1 \thalf mx^{l+2} e^{-mx}	N(x)\dx 
		\ \sim\ (-1)^l\frac{\hPhi(m)}{m^{l+2}}\,\frac{(l+2)!}2,
\ens
and
\eqs
  \lefteqn{\int_0^1 e^{m\log(1-x)}\{(-x)^l - \{\log(1-x)\}^l/(1-x)\}\,N(x)\dx}\\
	  &&\ \sim\ -\int_0^1 (-x)^l\thalf(l+2)xe^{-mx} N(x)\dx \\
		&&\ \sim\ \thalf (l+2)(-1)^{l+1}(l+1)!m^{-l-2}\hPhi(m),
\ens
it follows that $|D_l(m)| = o(m^{-l-2}\hPhi(m))$ also.
The lemma now follows by expressing the differences $\Phi_0\ul(m)-\hPhi\ul(m)$
using the quantities $(D_j(m),\,0\le j\le l)$.
\endpf
\vskip0.5cm

\par 
In fact, the measure $\tnu$ is an object in its own right: it is 
the L{\'e}vy measure of the {\it geometric\/} or {\it multiplicative\/}
subordinator $\widehat{S}_t=1-\exp(-S_t)$.  
In terms of $\widehat S$ the composition ${\cal C}_n$ is defined as a record of occupancy counts for
the gaps in the range of $\widehat S$
that are hit by at least one atom of a homogeneous Poisson process on $[0,1]$ with rate $n$.


\subsection{Law of large numbers}

For a law of large numbers, 
we begin by noting that
\eq\label{ADB-mean}
  \ex K_n \Eq \ex A_n(\infty) \Eq \int_0^\infty \ex\hPhi(ne^{-S_t})\,{\rm d}t
=\int_0^\infty \Phi(ne^{-s})U({\rm d}s)
\en
where $U$ is the potential measure of $S$
({\it i.e.}
$U[0,s]$ is the expected time $S$ stays below $s$).
Now, since~$\hPhi$ is slowly varying, it is plausible that
\[
  \ex\hPhi(ne^{-S_t}) \ \sim \hPhi(ne^{-\ex S_t})	\Eq \hPhi(ne^{-t}).
\]
This motivates the introduction of
\eq\label{Psi-tilde}
  \hPsi(n)\ :=\ \int_0^\infty \hPhi(ne^{-t})\,{\rm d}t
	  \Eq \int_0^n \hPhi(t)\,\frac{{\rm d}t}{t}\,,
\en
as an approximation to~$K_n\,$; under A1 it follows
that $\Psi(n)\gg \log n$.
Our aim in this section is to show that in fact $K_n \sim \hPsi(n)$ 
for large~$n$.

\begin{lemma}\label{L-0.4} 
Assumption {\rm A1} implies $\hPsi(m)\gg\hPhi(m)$.
\end{lemma}

\begin{proof} This is a standard consequence of the
 uniform convergence theorem \cite{BGT}, which states that 
slow variation implies 
$
  \hPhi(mu)/ \hPhi(m) \to\ 1,
$
uniformly in $u$ bounded away from $0$ and $\infty$. 
\endpf\end{proof}

\vskip0.5cm
It hence follows that, for any fixed~$T$, and as $\nti$, 
\eq\label{ratio-1} 
  \frac{A_n(T)}{\hPsi(n)}\ <\ \frac{T\hPhi(n)}{\hPsi(n)}\ \to\ 0\quad {\rm a.s.}
\en
and also that
\eq\label{ratio-2}
   \frac{\int_0^T \hPhi(ne^{-t})\,{\rm d}t}{\hPsi(n)}
	   \ < \ {T\hPhi(n)\over \hPsi(n)}\ \to\ 0. 
\en
The convergence in (\ref{ratio-1}) and (\ref{ratio-2}) also holds 
if the fixed time~$T$ is replaced by
an a.s.\ finite random time~$\tau$ which is measurable with respect
to $\s\{S_t,\,t\ge0\}$.

\par The next 
lemma explores the error caused by replacing $U$ with the Lebesgue measure
in~\eqref{ADB-mean}.

\begin{lemma}\label{mean1} For an arbitrary subordinator~$S$ with 
${\mathbb E}S_1^2<\infty$, 
\eq\label{AVG-mean-bnd}
  {\Psi}(n)\leq {\mathbb E}\,K_n\leq {\Psi}(n)+\lambda {\Phi}(n)
\en
for some constant $\lambda>0$. It follows that, 
under {\rm A1}, $\ex K_n/\Psi(n) \to 1$ as $\nti$. 
\end{lemma}
{\bf Proof.}
Let $X$ be a random variable with density $N_0(x)$, independent of $S$ and $\Pi_n$.
The closed range $X+{\cal R}$ of the delayed subordinator $(X+S_t, t\geq 0)$ 
is a stationary counterpart of ${\cal R}$. Ignoring the interval $]0,X[\,$,
the expected number of  occupied gaps produced by $X+S_t$ is precisely
$$
  {\mathbb E}K_{ne^{-X}}=\Psi(n)
$$
because, by stationarity, the potential measure of $S_t+X$ is Lebesgue measure.
Since $K_n$ is nondecreasing in $n$, we have
$$
  K_n\geq_d  K_{ne^{-X}},
$$
and thence  ${\mathbb E}K_n\geq \Psi(n)\,.$
This inequality can also be argued analytically, by using $U[0,s]\geq s$ and 
the monotonicity of $\Phi$.

\par Let $\tau$ be the passage time for $S$ through $X$. Because $S$ has 
a nontrivial overshoot over $X$ (while $X+S$ has none, since $X+S_0=X$),
the number of occupied gaps within $[S_\tau,\infty]$ produced  by $S$ is stochastically smaller than 
$K_{ne^{-X}}$ (produced by $X+S$). It follows that
$$
  K_n-{\cal K}_n(\tau)<_d K_{ne^{-X}}\,;
$$
passing to expectations, we obtain 
$$
  {\mathbb E}K_n\leq {\mathbb E} K_{ne^{-X}}+{\mathbb E}{\cal K}_n(\tau)
	=\Psi (n)+\ex\int_0^\tau \Phi(ne^{-S_t}){\rm d} t
   \leq \Psi(n)+\Phi(n) {\mathbb E}\tau \,.
$$
Now, ${\mathbb E}\tau<\infty$ follows from ${\mathbb E}X<\infty$, which 
is implied by the assumption ${\mathbb E}S_1^2<\infty$.
This completes the proof of~\eqref{AVG-mean-bnd}, with $\l = \ex\t$.
The convergence of $\ex K_n/\Psi(n)$ is now a consequence of Lemma \ref{L-0.4}.
\endpf
\vskip0.5cm
\bigskip

\par We now show that~$K_n$ is
close to~$A_n(\infty)$.

\begin{lemma}\label{KAlink}
Under Assumption {\rm A1}
\[
  \ex(K_n - A_n(\infty))^2 \ \sim\ \hPsi(n).
\]
Furthermore, considering the whole path of $\KK_n - A_n$, we have
\[
	\lnti\pr\Bigl[\sup_{0 \le T \le \infty}|\KK_n(T) - A_n(T)| > b_n\Bigr]
	  \Eq 0,
\]
for all sequences~$b_n$ such that $b_n^{-2}\,\hPsi(n) \to 0$.
\end{lemma}

\begin{proof}
The difference ${\cal K}_n(T)-A_n(T)$ is a square-integrable
martingale of locally bounded variation with respect to the 
filtration~$\ff_{T,n}$, and has all jumps of size~$1$; this yields 
the formula \cite[Section 15.2]{LastBrandt}
\eq\label{K-A}
  {\mathbb E}(K_n-A_n(\infty))^2 \Eq {\mathbb E}A_n(\infty) 
	\ \sim\ \hPsi(n),
\en	
proving the first part.  The second follows from Kolmogorov's
inequality, which gives
\eq\label{K-A-path}
  \pr\Bigl[\sup_{0 \le T \le \infty}|\KK_n(T) - A_n(T)| > b_n\Bigr]
	  \Le b_n^{-2}\,\ex A_n(\infty).\quad\Box
\en
\end{proof}

\par Next is the law of large numbers for the compensator.

\begin{proposition}\label{P-0.5}
Under Assumption  {\rm A1}, as $n\to\infty$,
$
 A_n(\infty)/ \hPsi(n) \to\ 1
$
almost surely and in the mean.
\end{proposition}

\begin{proof}Fix $\varepsilon$, and define
\[
  \t \ :=\ \sup\{t\ge0:\,t^{-1}|S_t-t| > \e\},
\]	
finite almost surely. By the monotonicity of~$\hPhi$,
\eqs
  \lefteqn{(1+\e)^{-1}\int_{(1+\e)\t}^\infty \hPhi(ne^{-s})\,{\rm d}s \Eq
 \int_\tau^\infty \hPhi(ne^{-(1+\varepsilon)t}) \,{\rm d}t\ <\
 \int_\tau^\infty \hPhi(ne^{-S_t}) \,{\rm d}t}\\
  &<& \int_\tau^\infty \hPhi(ne^{-(1-\varepsilon)t}) \,{\rm d}t
	\Eq (1-\e)^{-1}\int_{(1-\e)\t}^\infty \hPhi(ne^{-s})\,{\rm d}s.
	\phantom{HHHHHHH}
\ens
Dividing by $\hPsi(n)$ and using \Ref{ratio-1} and~\Ref{ratio-2},
 we make the sandwich
$$
  {1\over 1+\varepsilon}\ <\ \liminf_{n\to\infty} \frac{A_n(\infty)}{\hPsi(n)}
	\Le \limsup_{n\to\infty} \frac{A_n(\infty)}{\hPsi(n)} \ <\ 
	{1\over 1-\varepsilon}\quad {\rm a.s.},
$$
and now let $\varepsilon \to 0$ to obtain almost sure convergence.
Convergence in the mean then follows from $\ex K_n=\ex A_n(\infty)$ 
and Lemma \ref{mean1}, together with the fact that $A_n(\infty)/ \hPsi(n) \ge 0$
a.s. 
\endpf\end{proof}
\vskip0.5cm

\par Finally, we have all ingredients to establish
the law of large numbers for~$K_n$.

\begin{theorem}\label{lln} 
As $n\to\infty$, ${K}_n/\hPsi(n)\to 1$  almost surely and in the mean.
\end{theorem}

\begin{proof} 
It follows from Lemma~\ref{KAlink} that
$$
  {\mathbb E}\left({K_n\over \hPsi(n)}-{A_n(\infty)\over \hPsi(n)}\right)^2 
	\Eq {1\over \hPsi(n)}\,{\mathbb E}\left({A_n(\infty)\over \hPsi(n)}\right)
	\ \sim\ {1\over \hPsi(n)}\,,
$$
and convergence in the mean follows from Proposition~\ref{P-0.5}. 
Then, because $\hPsi(n)$ is increasing, continuous and unbounded, 
we can select $n_j$ to satisfy $\hPsi(n_j)=j^{2}$. Then from 
$$
  \sum_{j=1}^\infty \left({K_{n_j}\over \hPsi(n_j)}-{A_{n_j}(\infty)\over \hPsi(n_j)}
	\right)^2<\infty
$$
we conclude, in a standard way, that 
$K_{n_j}/\hPsi(n_j)\to 1$ a.s. along the subsequence. Now, because both $K_n$ and $\hPsi(n)$ are increasing
in $n$, the inequalities
$$ { K_{n_j}\over\hPsi(n_{j+1})}\leq { K_{n}\over\hPsi(n)}\leq { K_{n_{j+1}}\over\hPsi(n_j)}\,$$
hold for $n_j\leq n\leq n_{j+1}$. The convergence almost surely follows from these relations and the trivial fact 
that $\hPsi(n_j)/\hPsi(n_{j+1})\to 1$.
\endpf\end{proof}

\vskip0.5cm
\par Along the same lines, $K_n\gg {\cal K}_n(T)$ for each fixed $T$. This property can be shown
to be characteristic  of the slow variation case.

\subsection{The variance}\label{variance}

The aim of this section is to derive asymptotics of the variance of $K_n$ and $A_n(\infty)$.
We start with the explicit formula 
\eq\label{var-expl}
  {\rm Var}\, A_n(\infty)=2\,\,\int_0^\infty \Phi(ne^{-s})U({\rm d}s)\int_0^\infty \Phi(ne^{-s-v})
    {\rm d}_v\{U[0,v]-U[s,s+v]\},
\en
where ${\rm d}_v$ indicates the active  variable of integration.
The formula is derived by using the representation of the path past $t$ as 
$(S_{u}, u\geq t)=_d (S_t+S_{u}', u\geq 0)$ with~$S'$ an independent copy of $S$, 
 and using a familiar symmetrisation trick for squared 
integrals:
\begin{eqnarray*}
  {\rm Var \,}A_n(\infty) &=&
     \ex \left(\int_0^\infty \Phi(ne^{-S_t}){\rm d}t\right)^2
		     -\left(\int_0^\infty \Phi(ne^{-s})\,U({\rm d}s)\right)^2\\
  &=& 2\, \ex \left(\int_0^\infty \Phi(ne^{-S_t})\,{\rm d}t
    \int_t^\infty \Phi(ne^{-S_{u}})\,{\rm d}u\right)\\
	&&\qquad	- 2\,\int_0^\infty \Phi(ne^{-s})\,U({\rm d}s)\int_s^\infty 
	  \Phi(ne^{-v})\, U({\rm d}v)\\
  &=& 2\, \ex \left(\int_0^\infty \Phi(ne^{-S_t})\,{\rm d}t
    \int_0^\infty \Phi(ne^{-S_t-S_{u}'})\,{\rm d}u\right)\\
	&&\qquad	- 2\,\int_0^\infty \Phi(ne^{-s})\,U({\rm d}s)
    \int_0^\infty \Phi(ne^{-s-v})\,U(s+{\rm d}v)\\
  &=&  2\, \int_0^\infty \Phi(ne^{-s})\,U({\rm d}s) \int_0^\infty \Phi(ne^{-s-v})\,
	  U({\rm d}v)\\
	&&\qquad	- 2\,\int_0^\infty \Phi(ne^{-s})\,U({\rm d}s)\int_0^\infty 
	  \Phi(ne^{-s-v})\, U(s+{\rm d}v)\\
  &=& 2\,\int_0^\infty \Phi(ne^{-s})U({\rm d}s)\int_0^\infty \Phi(ne^{-s-v})
    {\rm d}_v\{U[0,v]-U[s,s+v]\}.
\end{eqnarray*}
Next is a more informative asymptotic formula, very much in the spirit of the asymptotics 
for the expectation $\ex K_n\sim \Psi(n)$ derived before.  To state it, we first define
\eq\label{psi2-def}
  \Psi_2(n) \ :=\ \int_0^\infty \hPhi^2(ne^{-s})\ds 
     \Eq \int_0^n {\Phi^2(s)}\,{{\rm d}s\over s}\,.
\en     

\begin{lemma}\label{var1} For a subordinator such that~$\hPhi$ is slowly varying at infinity, 
\eq\label{varA}
  {\rm Var}\,A_n(\infty)\ \sim\ \s^2\Psi_2(n)\,.
\en
If also {\rm A1} holds, then 
the same asymptotics hold for $K_n$;
$$ 
  {\rm Var}\,K_n\ \sim\ {\rm Var}\, A_n(\infty)\ \sim\ \s^2\Psi_2(n)\,.
$$
\end{lemma}

\begin{proof} 
Renewal theory  tells us that
$$
  0\ \leq\ U[0,s]-s\ \to\ \sigma^2/2\,,
$$
where the constant appears as the mean value of the delay variable~$X$: 
$$
 \ex X = \int_0^\infty x N_0(x)\,{\rm d}x
   ={1\over 2}\int_0^\infty x^2\,\nu_0({\rm d}x)={1\over 2}{\rm Var}\,S_1=
   {\sigma^2\over 2}\,
$$
(the third equality follows from (\ref{LK})).
This general fact holds for arbitrary square-integrable 
subordinators, and it  follows easily from  
the compound Poisson case treated in \cite{Feller}. 

\par We can therefore, for $\e$ given, select $s_0$ and $v_0$ so large that,
for $s>s_0$ and $v>v_0$,
\eq\label{Uest}
  |U[0,v_0]-U[s,s+v_0]-\sigma^2/2|<\e\,;~~~~~
   |U[v_0,v]-U[s+v_0,s+v]|<\e\,.
\en
Now, writing $G_s(v) = U[0,v] - U[s,s+v]$, we have 
\[
  G_s(0) \Eq 0,\quad \|G_s\| \ :=\ \sup_{v>0}|G_s(v)| \Le \s^2 \quad\mbox{and}\quad
    |G_s(v) - \s^2/2| \ <\ \e
\]
for $s>s_0$ and $v>v_0$.  Hence, by partial integration, the inner integral is at most
$\|G_s\|\hPhi(ne^{-s})$, so that truncating the external integral in 
(\ref{var-expl}) at the lower bound $s_0$ yields an error of
at most $ \lambda s_0\Phi^2(n)$, which is negligible when compared with the 
claimed asymptotics. Similarly, truncating the external integral at an upper bound 
$\log n-\psi$ yields an error of at most
\[
  2\s^2\int_{\log n-\ps}^\infty \hPhi^2(ne^{-s})\,U({\rm d}s)
    \Le 2e^{2\ps}\s^2 \Blb\int_0^\infty e^{-2v}\dv + U^*\Brb \Eq 2e^{2\ps}\s^2(1+U^*),
\]
using~\eqref{Phi-bnd-1}, where $U^* := \sup_{s>0} |U(s)-s|$.  Then    
truncating the internal integral at the upper bound $v_0$ results in an error  estimated as
\eq\label{epsup}
  2\,\int_{s_0}^\infty \Phi(ne^{-s})\,U({\rm d}s)\int_{v_0}^\infty \Phi(ne^{-s-v})\,
    {\rm d}\{U[v_0,v]-U[s+v_0,s+v]\}< 2\e 
    \,\int_{s_0}^\infty \Phi^2(ne^{-s})\,U({\rm d}s), 
\en
while, as above,
\[
  \Blm \int_0^\infty \hPhi^2(ne^{-s})(U({\rm d}s) - {\rm d}s)\Brm \Le 2U^*\hPhi^2(n),
\]
so that  
\eq\label{Ud}
\int_0^\infty \Phi^2(ne^{-s})\,U({\rm d}s)\sim \int_0^n \Phi^2(s)\,{{\rm d}s\over s}.
\en
\par 
Thus we are reduced to evaluating
$$
  J_{\e,n} :=2\,\int_{s_0}^{\log n-\psi_n} \Phi(ne^{-s})\,U({\rm d}s)\int_{0}^{v_0} \Phi(ne^{-s-v})\,
   {\rm d}\{U[0,v]-U[s,s+v]\},
$$
where we let $\psi=\psi_n\to\infty$ slowly enough that 
\[
  e^{2\ps_n} \ \ll\ \int_0^n \Phi^2(ne^{-s})\,U({\rm d}s).
\]
 But in the range $v<v_0$, $s<\log n-\psi_n$ we have $e^{-v}$  bounded from $0$ and $\infty$, and
$ne^{-s}>e^{\psi_n}\to\infty $;
hence, by the uniform convergence theorem for slowly varying functions \cite{BGT},
$$
  {\Phi(ne^{-s-v})\over \Phi(ne^{-s})}\to 1\,,~~~~~~{\rm as~~}n\to\infty,
$$
uniformly in such $s$ and $v$. 
With this substitution and using (\ref{Uest}) and (\ref{Ud}) we obtain
$$
  J_{\e,n} \sim (1+O(\e))\sigma^2 \int_{s_0}^\infty \Phi^2(ne^{-s})\,U({\rm d}s)\,
    \sim (1+O(\e))\sigma^2 \int_0^n \Phi^2(s)\,{{\rm d}s\over s}\,.
$$ 
Hence, recalling (\ref{epsup}) and sending $\e\to\ 0$, the desired asymptotics follow. 

\par Noting that, under~A1, 
\eq\label{ADB:small-var}
  \Psi_2(n) \Eq \int_0^n \Phi^2(s)\, {{\rm d}s\over s}\ 
	\gg\ \int_0^n \Phi(s)\, {{\rm d}s\over s} \Eq \Psi(n),
\en
it follows from Lemma~\ref{KAlink} that the asymptotics of ${\rm Var} \, K_n$  
are implied by those of ${\rm Var}\, A_n(\infty)$.
\endpf\end{proof}

\vskip0.5cm
\noindent
{\bf Remarks.}
A general subordinator~$\tS$ with
\[
 \mm\ :=\ \ex \tS_1\ =\ \int_0^\infty x\,\tinu_0({\rm d}x)\quad\mbox{and}\quad
   \t^2 \ :=\ \var\tS_1 \Eq  \int_0^\infty x^2\,\tinu_0({\rm d}x)
\]
yields a subordinator~$S$ with $\ex S_1 = 1$, by scaling time so that $S_t = \tS_{t/\mm}$;
$S$ has $\hPhi = \mm^{-1}\tPhi$ and $\s^2 = \mm^{-1}\t^2$, and has the same quantities
$K_n$ and~$A_n(\infty)$ as~$\tS$.  
Hence, for~$\tS$, we have
\[
  \ex A_n(\infty) \ \sim\ \int_0^n s^{-1}\hPhi(s)\ds \Eq \mm^{-1}\int_0^n s^{-1}\tPhi(s)\ds
\]
and
\[
  \var A_n(\infty) \ \sim\ \s^2\int_0^n s^{-1}\hPhi^2(s)\ds \Eq \t^2\mm^{-3}\int_0^n s^{-1}\tPhi^2(s)\ds\,.
\]  
For the gamma subordinator with $\tinu_0({\rm d}x) = 
a x^{-1}e^{-\th x}\dx$,  we have $\tPhi_0(s) = a\log(1 + s/\th)$, $\mm = a/\th$
and $\t^2 = a/\th^2$. Hence
\[
  \var K_n \ \sim\ \var A_n(\infty) \ \sim\ \th\int_0^n s^{-1} \log^2(1 + s/\th)\ds
    \ \sim\ \frac{\th}3 \log^3n,
\]
agreeing with the asymptotics    
for gamma-like subordinators obtained in \cite{GPYII} by a method based on the Mellin transform.
 
\par In the compound Poisson case, 
the asymptotics of $\var A_n(\infty)$ and $\var K_n$ are different, because~$\Psi(n)$
is no longer of smaller order than~$\Psi_2(n)$ as in~\eqref{ADB:small-var}. Instead,
with the normalisation $\nu_0[0,\infty]=1$, so that $\lnti \tPhi(n)=1$,  we have 
$$
  {\rm Var}\, A_n(\infty)\ \sim\ {\t^2\over {\tt m}^3}\,\log n\,;~~~~~
    {\rm Var}\, K_n \ \sim\ {\t^2 - {\tt m}^2\over {\tt m}^3}\,\log n,
$$
(see \cite{Bernoulli}), so that~(\ref{varA}) is valid only for the compensator.

\par In the case of regular variation, $\Psi_2(n)$
still gives the correct order of growth for the variances of  both quantities,
but the coefficients are not as in (\ref{varA}); see 
\cite{GPYI} for details.

It is immediate from the definitions \eqref{Psi-tilde} and~\eqref{psi2-def}
that $\liminf_{\nti} \Psi_2(n)/\Psi(n) \ge \hPhi(M)$ for all $M>0$, and hence
that $\Psi_2(n) \gg \Psi(n)$ as $\nti$. In consequence, from \eqref{varA}
and Lemma~\ref{KAlink}, 
\ignore{
The asymptotic formulas for the moments 
(\ref{meanKas}) and (\ref{varA}) are instances of the same
  integral
$$\Psi_q(n):= \int_0^\infty \Phi(ne^{-s})^q\,{\rm d}s=\int_0^n {\Phi(s)^q\over s} \,\,{\rm d}s$$
which converges for all $q>0$. The order of growth of $\Psi_q(n)$, as $n\to\infty$, 
 strictly increases with $q$.
This observation, contrasting $q=1$ and $2$, is especially important 
}
the fluctuations of the process $\KK_n  - A_n$ are of smaller order than those of~$A_n$,
so that, when studying limit theorems for~$\KK_n(t)$, it is enough to consider~$A_n(t)$.
\ignore{
Indeed, from Lemma~\ref{KAlink}, the extreme fluctuations of the process
$\KK_n  - A_n$  are at most of order $O(\hPsi(n)^{1/2})$, with $\Psi=\Psi_1\,$,
which is of the smaller order than 
$(\var K_n)^{1/2}\asymp \Psi_2(n)^{1/2}$.
In consequence, the additional fluctuations of~$\KK_n$ around~$A_n$
are negligible in comparison, 
and, as far as limit theorems are concerned, we can use~$A_n$ in place of~$\KK_n$.
}

\section{The key assumption}\label{key}

The asymptotics of moments only required the monotonicity of $\Phi$ and the  property
of slow variation.
In order to progress to a finer description of the asymptotics
of~$K_n$, we need  
a further assumption in addition to A1.
To express it, 
we begin by associating with~$\hPhi$ the function
$$
  L(s)\ :=\ {\hPhi(s)\over s\hPhi'(s)}\,,
$$
so that $(sL(s))^{-1}=(\log\hPhi(s))'$ and
\eq\label{Phi-L}
  \hPhi(s) \Eq \hPhi(1)\exp\left\{\int_1^s {{\rm d}z\over zL(z)}\right\}.
\en
Assumption A1 forces $\lim_{s\to\infty}L(s)= \infty$, because
the last formula is just 
an instance of the Karamata representation for slowly varying functions \cite{BGT}. 
Constantly keep in mind that the faster $L$, the slower $\Phi$.

\par Our extra assumption on $\hPhi$ is expressed via $L$, and puts
a limit on the way in which it can vary locally: we assume that 
there exist $s_0 \ge 1$ and $k>0$ such that
\begin{equation}\label{keyassum}
   \phantom{}\qquad{\bf Assumption\  A2}\colon\qquad  
	  \left|{sL'(s)\over L(s)}\right|\ <\ {k\over \log s}
	\quad \mbox{for all}\quad s\ge s_0.\phantom{H}
\end{equation}

\noindent
Because the right side in (\ref{keyassum}) goes to zero with $s$, 
the function $L$ is itself slowly varying;
under A1 the latter property is equivalent to the slow variation of $s\Phi'(s)$.

\begin{lemma}\label{2phisder}
If {\rm A1} holds and $L$ is slowly varying, then for $L_0(s):=\Phi_0(s)/(s\Phi_0'(s))$
$$
  \Blm \frac{sL'(s)}{L(s)} - \frac{sL_0'(s)}{L_0(s)} \Brm = o(s^{-1}L^2(s));
$$
thus Assumption~{\rm A2} can equivalently be stated using $L_0$ in place of~$L$.
\end{lemma}

\begin{proof} 
Direct calculation shows that
\[
  \frac{sL'(s)}{L(s)} - \frac{sL_0'(s)}{L_0(s)} =
	  s\Blb\Bl \frac{\hPhi'(s)}{\hPhi(s)} - \frac{\Phi_0'(s)}{\Phi_0(s)} \Br
		- \Bl \frac{\hPhi''(s)}{\hPhi'(s)} - \frac{\Phi_0''(s)}{\Phi_0'(s)} \Br
    \Brb.
\]
Now we apply Lemma~\ref{2phis} to bound differences between the derivatives
of~$\Phi_0$ and those of~$\hPhi$, and~\eqref{Phi-deriv-bnd} to bound the
derivatives themselves, and we also note that $\hPhi'(s) = \hPhi(s)/sL(s)$.
The lemma follows.
\endpf\end{proof}
		
\medskip We note in passing that both functions $1/L_0$ and $1/L$ can be given various 
probabilistic interpretations.  
For instance, in the spirit of (\ref{LK}),
$$1/L_0(n) =\ex \exp\{-n(S_\tau-\xi)\}$$
where $\xi$ is an independent exponential level with rate $n$ and $\tau$ 
is the passage time across $\xi$,
so that $S_\tau-\xi$ is the overshoot at $\xi$, see \cite[Corollary 1 (ii)]{Winkel}.
The function $1/L$ determines a conditional rate for creating singleton blocks of ${\cal C}_n$, meaning that
 $1/L(s)$, with $s=ne^{-S_t}$, is the conditional probability that 
a jump of $S$ at time $t$ covers exactly one Poisson point given at least one point
is covered.

\par Although we regard A2 as a {\it local} condition, under circumstances 
it can restrict the {\it global} growth of $L$.
For suppose that~$L$ is eventually increasing. Then, introducing 
$h(m):=L(m)/(mL'(m))$, we have in the usual way
$$L(m)=L(1)\,\exp\left( \int_{1}^m {{\rm d}z\over z h(z)}
\right).$$
Now, because eventually $L'\geq 0$, Assumption A2  reads as 
$h(m)>k^{-1}\log m$, hence yielding a global bound $L(m)< \lambda \log^k m$. 
In the other direction, observe that, even if $L$ is not monotone,
 the inequality inverse to (\ref{keyassum}), 
$L(m)>\lambda \log^{k} m$ with some $k>1$, would  disagree with A1, 
because in this case $\Phi$ would be bounded.

\bigskip
\noindent
{\bf Remark.} Of course, when (\ref{keyassum}) holds for some $s_0$ and $k$, we can set 
$s_0=1$ by taking $k$ sufficiently large.
This will suffice for our purposes, but gives a poor idea of the growth of~$L$.

\bigskip

Assumption A2 implies that 
\begin{equation}\label{bounds}
  \left({\log y\over \log x}\right)^{-k}\ <\ {L(y)\over L(x)}\ <\ 
		\left({\log y\over \log x}\right)^{k}\,,~~~~~~~s_0<x<y\,.
\end{equation}
To see this, for $s_0<x<y$, observe that
$$
  {L(y)\over L(x)}\Eq\exp\left\{\int_x^y {{\rm d}z\over h(z)\,z}
\right\}
\ < \
     \exp\left\{\int_x^y {{\rm d}z\over z}\left({k\over \log z}\right)\right\}
     \Eq\left({\log y\over \log x}\right)^{k},
$$
and similarly that
$$
  {L(y)\over L(x)}\ >\
    \exp\left\{-\int_x^y {{\rm d}z\over z}\left({k\over \log z}\right)\right\}
    \Eq\left({\log y\over \log x}\right)^{-k}.
$$

\par
Assumption A2 is a kind of `second order' slow variation, in the sense that 
the function $\log \hPhi$ has a form of de Haan's property (see \cite[Section 3.0]{BGT}):
$$
  {\log \hPhi(ne^{-s}) -\log \hPhi(n)\over 1/L(n)}\ \to\ -s\,,~~~~~~(n\to\infty).
$$
However, A2 is stronger than just this, and  
offers a better control on the variability of~$\Phi$; in particular we 
have the following estimate for 
the remainder.

\begin{lemma}\label{keylemma} 
Under assumptions {\rm A1--A2}, we have
\begin{equation}
  \left| \log\left\{ {\hPhi(me^{-s})\over \hPhi(m)}\right\}+{s\over L(m)}\right|\
	  <\ {\varkappa s^2 \over L(m)\log m}\,,
\end{equation}
for all $m\ge1$ and $|s|<\half\log m$, where $\k = k2^k$.
\end{lemma}

\begin{proof}
By Taylor's formula with the remainder in Lagrange's form,
$$
  \log\hPhi(me^{-s}) \Eq \log \hPhi(m)-\,{s\over L(m)}
	  -\,{s^2\over 2}{m^*L'(m^*)\over L(m^*)^2},
$$
for some~$m^*$ such that 
$$
  m\wedge(me^{-s})\leq m^*\leq m\vee (me^{-s})\,.
$$
Since $|s|<(\log m)/2$, we have $m^{1/2}\leq m^*\leq m^{3/2}$, and hence 
$\half\log m \le \log m^*\leq \thrhalf\log m$; and also, from A2,
$$
  \Blm{m^*L'(m^*)\over L(m^*)}\Brm\ <\ {k\over \log m^*}\Le {2 k\over \log m}\,.
$$
On the other hand, for $m^*$ in this range, 
(\ref{bounds}) implies that $L(m^*)> 2^{-k} L(m)$.  Hence
$$
  {m^*L'(m^*)\over 2 L(m^*)^2}\ <\ {k2^k\over  L(m)\log m}\,,
$$
as required. \hfill$\Box$
\end{proof}

\bigskip
\nin{\bf Remark.} So, loosely speaking, we are dealing with functions~$L$ 
that grow slowly enough, satisfying $L(n^{1/2})\asymp L(n)$. 
Indeed, it can be shown that $L(n)\gg L(n^{1/2})$ implies the convergence
$$
  \int_1^\infty {{\rm d}s\over sL(s)}\ <\ \infty\,,
$$
in which case the L{\'e}vy measure is finite. 

\begin{corollary}\label{ADB-cross}
Under Assumptions {\rm A1--A2}, we have
\[
  \hPhi(n)\, e^{-3s/2L(n)} \Le \hPhi(ne^{-s}) \Le \hPhi(n)\, e^{-s/2L(n)}
\]
for  $0 \le s \le \tfrac1{2(\k\vee1)}\,\log n$.
\end{corollary}	

\begin{proof}
Immediate from the above.
\endpf\end{proof}

\begin{corollary}\label{ADB-Psi(n)-bds}
Under Assumptions {\rm A1--A2}, if $ L(n) < \tfrac1{4\kvi}\log n$, we have
\ignore{\[
  \hPsi(n) \Le 
	\begin{cases}
	1 + \hPhi(n)L(n)\left(2 + \tfrac{4\kvi}e\right) 
	   &\mbox{{\rm if} } L(n) < \tfrac1{4\kvi}\log n;\\
	1 + \hPhi(n)\log n\left(\tfrac1{2\kvi} + \tfrac1e\right) 
	   &\mbox{{\rm if} } L(n) \ge \tfrac1{4\kvi}\log n\,.
	\end{cases}
\]
}
\eqs
  &&\hPsi(n) \Le 1 + \hPhi(n)L(n)\left(2 + \tfrac{4\kvi}e\right);\\
  &&\tfrac13 (1-e^{-6})\hPhi^2(n)L(n) \Le \Psi_2(n) 
     \Le \half + \hPhi^2(n)L(n)\{1 + \tfrac{2\kvi}e\}.
\ens
If  $ L(n) \ge \tfrac1{4\kvi}\log n$, we have    
\eqs
  &&\hPsi(n) \Le 1 + \hPhi(n)\log n\left(\tfrac1{2\kvi} + \tfrac1e\right);\\
  &&\tfrac13 (1-e^{-6})\hPhi^2(n)\log n \Le \Psi_2(n) 
     \Le \half + \hPhi^2(n)\log n\{\tfrac1{2\kvi} + \tfrac1{e^2}\}.
\ens
\end{corollary}

\begin{proof}
Write $k_n := \tfrac1{2\kvi}\log n$.
From the upper bound in Corollary~\ref{ADB-cross}, we have
\[
  \int_0^{k_n} \hPhi(ne^{-t})\,\dt \Le \hPhi(n)\{2L(n) \wedge k_n\},
\]
and
\eqs
  \int_{k_n}^{\log n}	\hPhi(ne^{-t})\dt &\le& \hPhi(n)e^{-k_n/2L(n)}\log n	\\
	&=& 4\kvi\hPhi(n)L(n) \Blb \frac{k_n}{2L(n)} \exp\Bl-	\frac{k_n}{2L(n)}\Br\Brb;
\ens
furthermore, from \Ref{Phi-bnd-1},
\[
  \int_{\log n}^\infty \hPhi(ne^{-t})\dt 
	  \Le \int_{\log n}^\infty ne^{-t}\dt	 \Eq 1.
\]
The bounds for~$\hPsi(n)$ now follow from its definition, and
because $xe^{-x} \le e^{-1}$ for $x\ge0$.  The proof of the upper bounds
for $\Psi_2(n)$ is analogous.

\par For the lower bound on~$\Psi_2(n)$, integrate~$\hPhi^2(ne^{-t})$ from $0$
to~$\min\{k_n,2L(n)\}$, and then use the lower bound in Corollary~\ref{ADB-cross}.	
\endpf\end{proof}

\ignore{
This last corollary is important when studying the fluctuations of~$\KK_n(t)$,
because it typically allows us to replace~$\KK_n(t)$ by~$A_n(t)$.
Indeed, from Lemma~\ref{KAlink}, the extreme fluctuations of the process
$\KK_n  - A_n$  are at most of order $\{\hPsi(n)\}^{1/2}$,
and  the scale of fluctuations of~$A_n$ has order $\{\Psi_2(n)\}^{1/2}$,
from Lemma~\ref{var1}.
Now, from Corollary~\ref{ADB-Psi(n)-bds}, it follows that
\[
  \hPsi(n) \Eq O\Bl \hPhi(n)\{\log n \wedge L(n)\} \Br,
\]
whereas
\eq\label{fluct}
  \Psi_2(n) \ \asymp\ \hPhi(n)^2	\{\log n \wedge L(n)\}.
\en
In consequence, the additional fluctuations of~$\KK_n$ around~$A_n$
are negligible in comparison, because $\hPhi(n)\to\infty$, and, as
far as limit theorems are concerned, we can use~$A_n$ in place of~$\KK_n$.
}

\bigskip
For the rest of this paper both assumptions A1 and A2 will be taken for 
granted, even if not explicitly mentioned.

\section{The forward argument}\label{forward}
In this section, under a wide range of circumstances
in which $L(n) = O(\log n)$, we show that the quantity	
\eq\label{An-star}
  A_n^*(T) := \int_0^{T\wedge\log n} \hPhi(ne^{-t})
    \Bl 1 - \frac{S_t-t}{L(ne^{-t})} \Br \dt
\en
is an adequate approximation to~$A_n(T\wedge\t_n)$, and hence,
in view of Lemma~\ref{tau-truncation}, to~$A_n(T)$.  
This is a very attractive result, because the random process~$S$ appears
only linearly in $A_n^*(T)$, making it easier to determine
the approximate behaviour of $A_n(T)$  from knowledge of that of~$S$. 
The way that the approximation is proved is to show that the quantity
$\sup_{T\ge0}|A_n(T\wedge\t_n)-A_n^*(T)|$ is asymptotically smaller
than the scale of fluctuations of~$A_n$.
In view of Lemma~\ref{tau-truncation} and Corollary~\ref{ADB-Psi(n)-bds}, 
in order to achieve this
when $L(n) = O(\log n)$,
we need to prove that, with probability tending to~$1$,
\[
  \sup_{T\ge0}|A_n(T\wedge\t_n)-A_n^*(T)| = o(\hPhi(n)\sqLn).
\]
To this end, we define the centred process
$$Z_t := S_t - t$$ 
and restrict attention
as far as possible to 
realisations
of~$S$ for which the paths of~$Z$ 
are reasonably nice.  This we make precise as follows.
First, for any $T, \f, \ps > 0$, we define the events
\eqa
 B_0(n) &:=& \{\half\log n \le \t_n \le 2\log n\}\label{B0-def};\\
 B_1(T,\f) &:=& \Blb \sup_{0\le t\le T} (t\vee1)^{-1/2}|Z_t| \le \f\Brb;
  \label{B1-def}\\
 B_2(\ps) &:=& \Blb\sup_{t\ge\ps} 2t^{-1}|Z_t| \le 1 \Brb. \label{B2-def}
\ena
The paths of~$Z$ are well behaved if $B_1(T,\f)$ holds for~$\f$ not too
large and for large enough~$T$, and if $B_2(\ps)$ holds for~$\ps$ not
too large.  With reference to these desiderata,
we have the following lemma.

\begin{lemma}\label{B-probs}
For $T, \f, \ps > 0$, we have
\eqa
  \pr[B_0^c(n)] &\le& 8\s^2/\log n;\label{B0-bnd}\\
  \pr[B_1^c(T,\f)] &\le& 2\s^2 \f^{-2}\lceil \log_2T \rceil \label{B1-bnd}~~~~~~~(T>2);\\ 	
  \pr[B_2^c(\ps)] &\le& 32\s^2 \ps^{-1},  \label{B2-bnd}
\ena
and also $B_0(n) \supset B_1(2\log n,\half\log^\a n)$ for any $0 \le \a \le
\half$ and $n \ge 3$.
\end{lemma}

\begin{proof}
First, by Kolmogorov's inequality for the centred, independent increments
process~$Z$, we have
\[
  \pr[B_0^c(n)] \Le 
	  \pr\Bigl[\sup_{0 \le u \le 2\log n}|Z_u| > \half\log n\Bigr]
	  \Le 8\s^2/\log n.
\]
The remaining statements are proved by combining
Kolmogorov's inequality with geometric dissection,
in a rather standard fashion.	For the second inequality, we have
\eqs
  \pr[B_1^c(T,\f)] &\le& \sum_{r=1}^{\lceil \log_2T \rceil}
	  \pr\left[ \max_{0 \le t\le 2^r} (2^{r-1})^{-1/2}|Z_t| > \f \right]\\
	&\le& \sum_{r=1}^{\lceil \log_2T \rceil} \frac{2^r\s^2}{\f^2 2^{r-1}}
	  \ = \	2\s^2 \f^{-2}\lceil \log_2T \rceil.
\ens
For the third, we have
\eqs
  \pr[B_2^c(\ps)] &\le& \sum_{r\ge\lceil\log_2\ps\rceil} \pr\left[
	   \max_{0 \le t \le 2^r} 2|Z_t| > 2^{r-1} \right] \\
	&\le& \sum_{r\ge\lceil\log_2\ps\rceil} 16\s^2 2^{-r} \ \le\ 32\s^2\ps^{-1}.
\ens
\endpf
\end{proof}

The following corollary needs no proof.

\begin{corollary}\label{B-limits}
For any positive sequences $T_n,\f_n,\ps_n$ we have
\eqs
  {\rm (i)}&&\lim_{n\to\infty} \pr[B_1(T_n,\f_n)] = 1 \quad\mbox{if}\quad
	   \f_n^{-2}\log T_n \to 0;\\
  {\rm (ii)}&&\lim_{n\to\infty} \pr[B_2(\ps_n)] = 1 \quad\mbox{if}\quad
	   \ps_n \to \infty.
\ens
\end{corollary}

\medskip
For the further argument, we distinguish two cases, relating to the 
global pattern of growth of~$\hPhi(n)$, each of which needs separate 
treatment. The idea of the distinction can be seen from the following formula 
for the variance of the linearised compensator when $T<\log n$:
\eq\label{vari}
  {\rm Var}\,A_n^*(T)
	  = \sigma^2\,\int_0^T \{\Phi(ne^{-t})-\Phi(ne^{-T})\}^2 \,{\rm d}t\,,
\en
which is derived by writing (\ref{An-star}) for the centred $A_n^*$ 
as a stochastic integral:
$$
  A_n^*(T)-\ex A_n^*(T)
	  =  -\int_0^T \{\Phi(ne^{-t})-\Phi(ne^{-T})\}\,{\rm d}(S_t-t)\,,
$$
and using the independence of increments.
So, when $\Phi$ is a function like a power of logarithm, the difference
$\Phi(ne^{-t})-\Phi(ne^{-T})$ is of constant order over the 
whole time-range from~$0$ to $\log n$. On the other hand, if 
$\Phi$ grows fast enough, the first term will dominate, and the principal 
contribution to the integral will
come from times $t=o(\log n)$, as is also the case if~$\hPhi$ is regularly 
varying.

\subsection{Moderately growing $\hPhi$}
We begin with the boundary case, which includes the gamma-like subordinators~\cite{GPYII}, 
when $\hPhi(n)$ grows more or less like a power of~$\log n$.  Here,
all times~$t$ between $0$ and~$\log n$ contribute more or less
evenly to the fluctuations of~$A_n$. This case is defined by a
global condition on the function~$L$; that, for 
some $1 \le c_2 < \infty$ and for some~$m_0$, and with $c_1 := \{3\kvi\}^{-1}$,
\eq\label{L-cond1}
   \frac{c_1\log m}{6\log\log m} \Le L(m) \Le c_2\log m,\quad m \ge m_0.
\en
The next lemma is a preliminary to proving that, under these
circumstances, $A_n^*$ is a good approximation
to~$A_n$.  It enables us to truncate the integrals 
defining $A_n(T)$ and $A_n^*(T\wedge\t_n)$ close to~$\log n$, when
the paths of~$Z$ are nice enough.

\begin{lemma}\label{L4}
On the event $B_1(2\log n,\f_n)$, and for $0 < \ps_n \le \log n - \f_n\sqln$,
we have, for all $T>0$,
\eqa
  &&0 \Le \int_0^{\tnT} \hPhi(ne^{-S_t}) \dt - 
    \int_0^{\psnT} \hPhi(ne^{-S_t}) \dt \Le \h_n;
		\label{L4-E1}\\
  &&\int_{\psnT}^{T\wedge\log n} \hPhi(ne^{-t})
	   \Blm 1 - \frac{Z_t}{L(ne^{-t})}\Brm \dt \Le \h_n,
		 \label{L4-E2}
\ena
where		 
\[		
  \h_n = (\log n - \ps_n + \f_n\sqrt{2\log n})
	  \hPhi(n\exp\{-\ps_n + \f_n\sqpsn\}).
\]
\end{lemma}

\begin{proof}
On $B_1(2\log n,\f_n)$, we have
\[
  \log n - \f_n\sqrt{\log n} \Le \t_n \Le \log n + \f_n\sqrt{2\log n},
\]
and so $\ps_n \le \min\{\t_n,\log n\}$.
Hence \Ref{L4-E1} and~\Ref{L4-E2} are both zero 
if $T \le \ps_n$. The first part of the lemma then merely uses the fact that 
\[
  \ps_n \le \t_n \le \log n + \f_n\sqrt{2\log n},
\]
combined with the largest possible value of the integrand in this range.
For the second part, recall (\ref{gnt})
so that
\[
  \int_{\ps_n}^{\log n} \hPhi(ne^{-t})\,\dt \Le \hPhi(ne^{-\ps_n})(\log n - \ps_n)
\]
and from (\ref{B1-def})
\[\phantom{HHHHHH}
  \int_{\ps_n}^{\log n} \hPhi(ne^{-t})\frac{|Z_t|}{L(ne^{-t})} \,{\rm d}t
	  \Le \hPhi(ne^{-\ps_n})\f_n\sqrt{\log n}.\phantom{HHHHHH}\Box
\]
\end{proof}

\medskip
It follows from~A2, (\ref{L-cond1}) and the definition of $L$ 
that
\eq\label{Nhat-bnd2}
  \hPhi(ne^{-t}) \Eq \hPhi(n)\exp\Blb -\int_{ne^{-t}}^n \frac{dy}{yL(y)}\Brb 
  \Le \hPhi(n)\Blb 1 - \frac{t}{\log n}\Brb^{1/c_2}\phantom{HH}  
\en
for all $n$ and~$t$ such that $ne^{-t} \ge m_0$. Thus, taking
\eq\label{psi-def1}
  \ps_n = \log n - 2u_n\sqln,
\en
for $u_n \ge \f_n$,
the quantity~$\h_n$ in Lemma~\ref{L4} is, for all~$n$ large
enough, at most
\[
  \hPhi(n)\,4u_n\sqln\,\{3u_n/\sqln\}^{1/c_2}.
\]	  
This is in turn at most
\[
  12\sqrt{18\kvi}\,\hPhi(n)\sqLn\,\{\sqrt{\log\log n}
	\log^{-\b/2}n\}
\]	
if we take
$u_n = \log^{\b}n$ for $\b = 1/\{4(1+c_2)\}$. 

\begin{theorem}\label{T1}
Suppose that Assumptions {\rm A1--A2} and~\Refrm{L-cond1} hold, fix
$3\a =\b = \tfrac1{4(c_2+1)}$, and set $u_n = \log^\b n$,
 $\f_n = \log^\a n$.  Then, on $B_1(2\log n,\f_n)$,
we have
\[
  \sup_{T\ge0}\Blm A_n(T\wedge\t_n) - A_n^*(T)\Brm
 = \varepsilon(n)
\left(\hPhi(n)\sqLn\right),
\]
where $\lnti \e(n) = 0$ uniformly in
$c_2 < C < \infty$, for each $C>0$.  Furthermore,
\[
  \lnti \pr[B_1(2\log n,\f_n)] \Eq 0.
\]	
\end{theorem}

\begin{proof}
By the argument just completed, it is enough to examine the integrated difference
\[
  \int_0^{\psnT} \Blm \hPhi(ne^{-S_t})  
    -  \hPhi(ne^{-t}) \Blb 1 - \frac{Z_t}{L(ne^{-t})} \Brb \Brm \dt,
\]		
where $\ps_n = \log n - 2u_n\sqln$.

To this end,
we use Lemma~\ref{keylemma} with $ne^{-t}$ for~$m$ and $Z_t$ for~$s$.
On  $B_1(2\log n,\f_n)$, and since, for $0\le t\le\ps_n$, we have 
$ne^{-t} \ge \exp\{2u_n\sqln\}$, it
follows that 
\eqs
  |Z_t| &\le& \f_n\sqln \Eq \half\log\{\exp(2\f_n\sqln)\} \\
	&\le& \half\log\{\exp(2u_n\sqln)\} \Le \half\log\{ne^{-t}\},
\ens
so that the lemma can be applied. It then follows that
\eqa
\lefteqn{|\hPhi(ne^{-S_t}) - \hPhi(ne^{-t})\exp\{-Z_t/L(ne^{-t})\}|} \non\\
  &\le& \hPhi(ne^{-t}) \exp\{-Z_t/L(ne^{-t})\} X(n,t) \exp\{X(n,t)\}\,,
	\label{main-bnd1}
\ena
where
\eq\label{Xnt-def}
  X(n,t) \Eq \frac{\k Z_t^2}{L(ne^{-t})(\log n-t)}\,.
\en
It is then also immediate from $e^{-x}-1+x<e^{|x|}x^2/2$ that
\eqa
\lefteqn{\hPhi(ne^{-t})|\exp\{-Z_t/L(ne^{-t})\} - 1 + Z_t/L(ne^{-t})|} \non\\
  &&\le\ \half \hPhi(ne^{-t}) \exp\{|Z_t|/L(ne^{-t})\} \{Z_t/L(ne^{-t})\}^2\,.
  \label{main-bnd2}
\ena
Now, on  $B_1(2\log n,\f_n)$, and for  $0\le t\le\ps_n$, we have
\[
  X(n,t) \Le (6\k/c_1)(\f_n/u_n)^2\log\log n \Le \lambda_1,
\]	  
and also 
\[
  |Z_t|/L(ne^{-t}) \Le (6/c_1)(\f_n/u_n)\log\log n \Le \lambda_2;
\]	
more precisely, in this range of~$t$, by (\ref{L-cond1})
\[
  X(n,t) \Le \frac{6\k\f_n^2\log n\log\log n}{c_1(\log n - t)^2}
    \quad\mbox{and}\quad 
    \frac{|Z_t|}{L(ne^{-t})} \Le \frac{6\f_n\sqln\log\log n}{c_1(\log n - t)}\,.
\]
We also have the bound~\Ref{Nhat-bnd2} for~$\hPhi(ne^{-t})$.
Combining these, it follows that
\eqa
\lefteqn{\Blm \int_0^{\ps_n} \hPhi(ne^{-S_t}) \dt 
    -  \int_0^{\ps_n} \hPhi(ne^{-t})
    \Blb 1 - \frac{Z_t}{L(ne^{-t})} \Brb \dt\Brm} 
	  \non\\
  &\le& 
	 \int_0^{\ps_n} \hPhi(ne^{-t})
    \Blb e^{\lambda_1+\lambda_2}X(n,t) + \half e^{\lambda_2}(Z_t/L(ne^{-t}))^2 \Brb\dt
	  \phantom{HHHHHHH}\non\\
  &\le& \frac6{c_1^2}(\k c_1 e^{\lambda_1+\lambda_2} + 3e^{\lambda_2})
	  \f_n^2\hPhi(n)(\log\log n)^2\int_{2u_n/\sqln}^1 u^{-2+1/c_2}\du\label{2.30}\\
	&\le& \lambda \left(\hPhi(n)\sqLn\right) (\log\log n)^{5/2}\f_n^2/u_n \non\\
	&=& \lambda \left(\hPhi(n)\sqLn\right)\, \{(\log\log n)^{5/2}\log^{-\b/3}n\},
    \label{diff-bnd}
\ena 
for some $\lambda>0$. This completes the proof.  Note that~$c_2$ enters the
bound implicitly, in the value of~$\b$, and hence in $\lambda_1$ and~$\lambda_2$.
\hfill$\Box$
\end{proof}

\bigskip
\nin{\bf Remark.}\ The restrictions imposed by~\Refrm{L-cond1} can be
relaxed somewhat, to allow a little more freedom in  both lower and
upper bounds.  For instance, the same proof can be used under the
condition
\eq\label{L-cond2}
  (\log m)^{1 - \g(c_2(m))} \Le L(m) \Le c_2(m)\log m\quad
	  \mbox{for all}\quad m \ge m_0,
\en
for any increasing function~$c_2$ satisfying $c_2(m) = o(\log\log m/
\log\log\log m)$, where $\g(c) := 1/\{32(c+1)\}$. Suitable choices
of the parameters are now 
\[
  \b \Eq \b(n) \ :=\ 8\g(c_2(n)) \quad\mbox{and}\quad
	  \a \Eq \a(n) \ :=\ \b(n)/3;
\]
note that, with these definitions and with $\f_n = \log^{\a(n)}n$,  
we still have	$\pr[B_1^c(2\log n,\f_n)] \to 0$.	
This extra freedom enables the main transition, between the behaviour in the
case of moderately growing~$\hPhi(n)$ and that when~$\hPhi(n)$ grows 
either faster or more slowly, to be understood in greater detail.

\subsection{Fast growing $\hPhi$}
We  turn to the setting in which slowly varying~$\hPhi(n)$ grows faster
than any power
of~$\log n$. In this case most of the random fluctuation in~$A_n$
takes place at times of order~$L(n)$, where $L(n)$ goes to infinity 
(as required by A1)
but
slower than~$\log n$.  Our global condition determining this
r\'egime is
\eq\label{L-cond4}
\qquad 6L(n)\log L(n) \Le c_1\log n,
\en
where $c_1 = \{3\kvi\}^{-1}$ is as before.
Note that, if $L(n) \le {c_1\log n\over 6\log\log n}\,$, then (\ref{L-cond4})
is satisfied; the condition 
given in~\Ref{L-cond1} was chosen to match neatly, though
in view of the remark at the end of the previous section, this was
not really necessary. 
Here,
we first need a modification of 
Lemma~\ref{L4}, in order to be able
to truncate the integrals defining $A_n(T\wedge\t_n)$ and $A_n^*(T)$
as far as we need to.

\begin{lemma}\label{L4-new}
Suppose that $\ps_n$ is such that $6L(n)\log L(n) \le \ps_n \le c_1\log n$,
and that~$n$ is large enough to satisfy $L(n) \ge e^6$.
Then, on the event $B_2(\ps_n)$, we have
\eqa
  0 \Le \int_0^{\tnT} \hPhi(ne^{-S_t}) \dt - 
    \int_0^{\psnT} \hPhi(ne^{-S_t}) \dt 
	\ \le\ 4\hPhi(n)\{e^{-3} + c_1^{-1}\};\quad&&
		\label{L4-E1-new}\\
  \int_{\psnT}^{T\wedge\log n} \hPhi(ne^{-t})
	   \Blm 1 - \frac{Z_t}{L(ne^{-t})}\Brm \dt 
	\ \le\ \hPhi(n)\{2e^{-3} + \tfrac32c_1^{-1} + 6e^{-9} 2^k \}.\quad&&
		 \label{L4-E2-new}
\ena
\end{lemma}

\begin{proof}
On $B_2(\ps_n)$, we have $\tfrac23\log n \le \t_n \le 2\log n$, implying
immediately that $\ps_n \le \min\{\t_n,\log n\}$.  Hence, if $T\le\ps_n$,
both of the quantities to be bounded in the lemma are zero.  Note also, in
preparation, that for $l\ge e^6$ and for any $x \ge 6l\log l$, we have
\eq\label{xexl}
  xe^{-x/4l} \Le 6l\log l\exp\{-3\log l/2\} \Eq 6l^{-1/2}\log l \Le 2.
\en

For the bound \Ref{L4-E1-new}, since $S_t \ge t/2$ for $t\ge\ps_n$ 
on $B_2(\ps_n)$ and since $c_1 \le 1/2\kvi$, we can apply 
Corollary~\ref{ADB-cross} to give
\eqs
  \int_{\ps_n}^{c_1\log n} \hPhi(ne^{-S_t})\dt 
	  &\le&	\hPhi(n)\int_{\ps_n}^{c_1\log n} e^{-t/4L(n)}\dt \\
	&\le& 4L(n)\hPhi(n)e^{-\ps_n/4L(n)} \Le 4L(n)^{-1/2}\hPhi(n),
\ens
by the definition of~$\ps_n$.  Then we also have
\[
  \int_{c_1\log n}^{2\log n} \hPhi(ne^{-S_t})\dt
	  \Le 2\log n\, \hPhi(n) \exp\{-c_1\log n/4L(n)\} \Le 4c_1^{-1}\hPhi(n),
\]
this last by~\Ref{xexl}.

\par  The argument for~\Ref{L4-E2-new}	is very similar. First, bounding
$\int_{\ps_n}^{\log n} \hPhi(ne^{-t})\dt$, it follows from 
Corollary~\ref{ADB-cross} that
\[
  \int_{\ps_n}^{c_1\log n} \hPhi(ne^{-t})\dt 
	  \Le 2L(n)\hPhi(n)e^{-\ps_n/4L(n)} \Le 2L(n)^{-1/2}\hPhi(n),
\]
and then that
\[
  \int_{c_1\log n}^{\log n} \hPhi(ne^{-t})\dt
	  \Le \log n \,\hPhi(n) \exp\{-c_1\log n/2L(n)\} \Le c_1^{-1}\hPhi(n).
\]
For the remaining term, we first have
\eqa
  \int_{\ps_n}^{c_1\log n} \frac{\hPhi(ne^{-t})|Z_t|}{L(ne^{-t})}\dt
	  &\le& \frac{\hPhi(n)}{L(n^{1-c_1})}	\int_{\ps_n}^{c_1\log n} 
	  \half te^{-t/2L(n)}\dt\non\\
	&\le& \frac{2L(n)^2\hPhi(n)}{L(n^{1-c_1})} \int_{\ps_n/2L(n)}^\infty
	  ue^{-u}\du.\label{ADB-237}
\ena
Now, for any $y\ge 1/2$,
\[
  \int_y^\infty ue^{-u}\du \Eq \int_0^\infty (y+v)e^{-y-v}\dv \Eq
	  e^{-y}(y+1) \Le 3ye^{-y},
\]
so that~\Ref{ADB-237} can be bounded, using \Ref{bounds} and~\Ref{xexl}, by
\[
  \frac{3L(n)\hPhi(n)}{L(n^{1-c_1})}\,\ps_n e^{-\ps_n/2L(n)}
    \Le \frac{6L(n)\hPhi(n)}{L(n^{1/2})}\,e^{-\ps_n/4L(n)}
		\Le \frac{6\cdot 2^k\hPhi(n)}{\{L(n)\}^{3/2}}.
\]
Finally, using (\ref{gnt})
and 
\hbox{$c_1 < 1/2\kvi$,} we have
\eqs
  \int_{c_1\log n}^{\log n}	\frac{\hPhi(ne^{-S_t})|Z_t|}{L(ne^{-t})}\dt
	  &\le& \half\log n\,\hPhi(ne^{-c_1\log n})\\
	&\le& \half\log n\,\hPhi(n)\,e^{-c_1\log n/2L(n)} \Le \tfrac1{2c_1}\hPhi(n),
\ens
again using \Ref{xexl}.
This completes the proof.\hfill$\Box$
\end{proof}

\medskip
\begin{theorem}\label{VFGPhi}
Under Assumptions {\rm A1--A2} and~\Refrm{L-cond4}, set $\ps_n = 6L(n)\log L(n)$ 
and $\f_n = L(n)^{1/6}$.  Then, on the event $B_1(\ps_n,\f_n) \cap B_2(\ps_n)$,
and if $L(n) \ge e^6$, we have
\[
  \sup_{T\ge0}|A_n(T\wedge\t_n) - A_n^*(T)| 
	  \Le \e(L(n))\left(\hPhi(n)\sqLn\right)\,,
\]
where $\lim_{m\to\infty}\e(m) = 0$.
Furthermore,
\[
  \lnti \pr[B_1^c(\ps_n,\f_n)] \Eq \lnti \pr[B_2^c(\ps_n)] \Eq 0.
\]	
\end{theorem}

\begin{proof}
As before, on $B_2(\ps_n)$, we have $\tfrac23\log n \le \t_n \le 2\log n$, implying
immediately that $\ps_n \le \min\{\t_n,\log n\}$.	
By Lemma~\ref{L4-new}, it is enough to bound  the difference
\[
  \int_0^{\ps_n} \Blm \hPhi(ne^{-S_t})  
    -  \hPhi(ne^{-t}) \Blb 1 - \frac{Z_t}{L(ne^{-t})} \Brb \Brm \dt.
\]		
By~\Ref{bounds}, we can use the inequality $L(ne^{-t}) \ge 2^{-k}L(n)$
for $0\le t\le c_1\log n$.  Hence, on the event $B_1(\ps_n,\f_n)$,
and noting that $\log n - \ps_n \ge \half L(n)$ because of~\Ref{L-cond4},
we can bound the quantities $X(n,t)$ and~$\{Z_t/L(ne^{-t})\}^2$ appearing
in the proof of Theorem~\ref{T1} by
\[
  2^{2k}\f_n^2\ps_n L(n)^{-2} \Eq 6\cdot2^{2k}L(n)^{-2/3}\log L(n) 
	  \Le 36\,e^{-4}\,2^{2k}\,,
\]
in the range $t \le \ps_n$;  thus they are both
uniformly bounded in~$n$, and asymptotically small as $n\to\infty$.	
Hence, using \Ref{main-bnd1} and~\Ref{main-bnd2}, it follows that
\eqs
  \lefteqn{\int_0^{\ps_n} \Blm \hPhi(ne^{-S_t})  
    -  \hPhi(ne^{-t}) \Blb 1 - \frac{Z_t}{L(ne^{-t})} \Brb \Brm \dt} \\
  &&\Le \lambda L(n)^{-2/3}\log L(n) \int_0^{\ps_n}\hPhi(ne^{-t})\dt 
\ens
for some $\lambda<\infty$.	But now, from Corollary~\ref{ADB-cross}, it follows
that
\eqs
  L(n)^{-2/3}\log L(n) \int_0^{\ps_n}\hPhi(ne^{-t})\dt &\le& 
	  2L(n)^{1/3}\hPhi(n)\log L(n) \\
	&=& 2L(n)^{-1/6}\log L(n)\left(\hPhi(n)\sqLn\right)\,,
\ens
proving the main assertion.  The last statement follows from
Corollary~\ref{B-limits}.\hfill$\Box$
\end{proof}

\section{The backward argument}\label{backward}
We now turn to the case of functions~$\hPhi(n)$ that grow more slowly than
any power of~$\log n$.  Here, the argument required and the approximations
obtained are of rather different character to those of the previous
section.  In particular, we make use of properties of the L\'evy
process when looking backwards in time.  Our setting is defined by
requiring that $\lnti\hPhi(n) =\infty$, but that~$L$ satisfies the 
following global condition: 
\eq\label{L-cond6}
  L(m) = c_2(m)\log m, \quad\mbox{where}\quad 
	  \lim_{m\to\infty}c_2(m)\Eq\infty.
\en
To agree with A1, $c_2(m)$ must grow slowly enough, meaning that the 
integral in ~\Ref{Phi-L},
$$
  \int_2^\infty \frac{{\rm d}m}{c_2(m)m\log m}\,,
$$
must diverge, a condition which excludes functions 
like $c_2(m)=\log^\e m$ for any $\e>0$.
One can think of $c_2(m)
= \log\log m$ for $m \ge m_0$, as one possible example, in which case 
$\hPhi(n) \asymp \log\log n$. 
Here, we no longer have Lemma~\ref{L4} to help us.
However, the argument of Theorem~\ref{T1} is still good,
if we restrict to taking the supremum over $0\le T\le (1-\d_n)\log n$,
for some $\d_n \to 0$ sufficiently slowly, and this gives us the
following approximation of $A_n$ by~$A_n^*$.  

\begin{lemma}\label{T1-new}
Take $\a = 1/8$, $\f_n = \log^\a n$ and $\d_n=\log^{-1/8}n$.	
Then, on the event $B_1(2\log n,\f_n)$, it follows that
\[
  \sup_{0\le T\le (1-\d_n)\log n} |A_n(T\wedge\t_n) - A_n^*(T)|
	  \Le \lambda\,\frac1{\log^{1/8}n}\,\frac{\hPhi(n)\sqln}{c_2^*(n^{\d_n})}\,,
\]
for some $\lambda > 0$ and 
$c_2^*(m) = \inf_{r\ge m} c_2(r)$.
\end{lemma}

\begin{proof}
We argue as for Theorem~\ref{T1}, now with $\ps_n = (1-\d_n)\log n$,
noting that, for $t \le \ps_n$,
\[
  |X(n,t)|	\Le \frac{\k\f_n^2\log n}{(\log n-t)^2 c_2^*(n^{\d_n})}
	  \Le \frac{\k}{\f_n^2 c_2^*(n^{\d_n})}\,,
\]
since $\d_n > \f_n^2/\sqln$, and that
\[
  \frac{|Z_t|}{L(ne^{-t})} \Le \frac{\f_n\sqln}{(\log n - t) c_2^*(n^{\d_n})}
	  \Le \frac1{\f_n c_2^*(n^{\d_n})}\,,
\]					
both of which are small in~$n$. Then, arguing as for~\Ref{2.30}, and
using the crude bound $\hPhi(ne^{-t}) \le \hPhi(n)$, we have
\eqs
  |A_n(T\wedge\t_n) - A_n^*(T)| &\le& \lambda\,\frac{\hPhi(n)\f_n^2}{c_2^*(n^{\d_n})}
	  \int_{\d_n}^1 u^{-2}\du \\
	&\le& \lambda\,\frac{\hPhi(n)\f_n^2}{\d_n c_2^*(n^{\d_n})} 
	  \Le \lambda\,\frac1{\log^{1/8}n}\,\frac{\hPhi(n)\sqln}{c_2^*(n^{\d_n})}\,,
\ens
for $T\le (1-\d_n)\log n$, as required.\hfill$\Box$
\end{proof}		

\bigskip
To see that differences of this order are relatively small, we now make
some variance calculations, for which 
we introduce the notation 
\eq\label{W-def}
  -W(v) \ :=\ \int_0^{v\log n} g(n,t) (S_t-t)\dt ,
\en
where 
\eq\label{gnt}
  g(n,t)\ :=\ \hPhi(ne^{-t})/L(ne^{-t}) \Eq ne^{-t}\Phi'(ne^{-t})
  \Eq -\frac{{\rm d}}{{\rm d}t}\,\left\{\Phi(ne^{-t})\right\}\,.
\en
It thus follows that 
\eq\label{W-Astar}
  W(v) = A_n^*(v\log n) - \ex\{A_n^*(v\log n)\} 
\en	
whenever $vL(n) \le \log n$.

\begin{lemma}\label{Var-2}
For $0\le v \le 1$ and~$n$ large enough, we have
\[
  \var A_n^*(v\log n) \Ge \lambda(v\wedge\half)^3\hPhi(n)^2\log n/c_2(n)^2;
\]
for $0 < \d < 1$ and for $0 \le v \le (1-\d)$, we have
\[
  \var A_n^*(v\log n) \Le \lambda\hPhi(n)^2\log n\, \frac{\log(1/\d)}{c_2^*(n^\d)}\,,
\]
where $c_2^*(m) = \inf_{r\ge m} c_2(r)$.
\end{lemma}

\begin{proof}	
From Lemma~\ref{keylemma}, it follows that, for $0\le v \le \half$
and $t = v\log n$, 
\[
  \hPhi(ne^{-t}) \ge \hPhi(n)\, e^{-t/L(n)} \exp\{-\k v^2/c_2(n)\}\,,
\]
and, from~\eqref{bounds}, that $L(n) / L(ne^{-t}) \ge 2^{-k}$,	
implying that, for~$n$ so large that $\k/\{4c_2(n)\} \le 1$, we have
\[
  g(n,t) \ge  2^{-k}e^{-1}\hPhi(n)\, e^{-t/L(n)}/L(n)\,,
\]
where 
$g(n,t)$ is as in (\ref{gnt}).
%
Now $Z_t = S_t-t$ has independent increments with zero means, and $\var Z_t = 
\s^2t$.  Hence, for any~$0\le v\le 1/2$, recalling~\Ref{W-Astar}, we have
\eqs
  \var A_n^*(v\log n) &=& 2\int_0^{v\log n} \int_0^t  g(n,t)g(n,u) \s^2 u \du\dt\\
  &\ge& 2^{1-2k}e^{-2}\hPhi(n)^2 \int_0^{v/c_2(n)} \int_0^w  e^{-w-z} \s^2 zL(n) \dz\dw\\
	&\ge& 2^{1-2k}e^{-2}\s^2\hPhi(n)^2 c_2(n)\log n\cdot \tfrac16 (v/c_2(n))^3
	  e^{-1/c_2(n)}\\
	&\ge& 2^{1-2k}\s^2\hPhi(n)^2 c_2(n)\log n\cdot \tfrac16 (v/c_2(n))^3
	  e^{-3}\,,
\ens
for all~$n$ large enough.
This proves the first inequality, since this lower bound with
$v=1/2$ is a lower bound for larger~$v$ also.

For the second part, we recall~\Ref{vari}:
\[
   \var A_n^*(v\log n) 
	   \Eq \s^2\int_0^{v\log n} \{\hPhi(ne^{-t}) - \hPhi(n^{1-v})\}^2 \dt,
\]
whenever $0\le v\le 1$.  Now, from the representation~\Ref{Phi-L}, it
follows that, for $0 < T \le (1-\d)\log n$,
\eqs
  \lefteqn{\int_0^T \{\hPhi(ne^{-t}) - \hPhi(ne^{-T})\}^2 \s^2\dt}\\
	&&\le\ \int_0^T \hPhi^2(ne^{-t})\Blb 1 - 
	  \exp\Bl - \frac1{c_2^*(n^\d)}\int_{ne^{-T}}^{ne^{-t}} \frac{\dy}{y\log y}\Br
		\Brb^2\dt \\
	&&=\ \int_0^T \hPhi^2(ne^{-t})\Blb 1 - 
	  \Bl\frac{1 - T/\log n}{1 - t/\log n}\Br^{1/c_2^*(n^\d)}\Brb^2\dt \\
	&&\le\ T\,\hPhi^2(n)\{1 - \d^{1/c_2^*(n^\d)}\}\,,
\ens
and the second part is proved.	
\hfill$\Box$
\end{proof}

\bigskip
\nin In particular, the lower bound shows that the standard deviation 
of~$A_n^*(T)$ is at least
as big as a constant times $\hPhi(n)\sqln/c_2(n)$ for $T \ge \half\log n$.
By comparison, the differences in Lemma~\ref{T1-new} are typically much
smaller, because of the factor $\log^{-1/8}n$;
recall that $c_2(n)$ grows rather slowly with~$n$, and certainly not as
fast as a power of~$\log n$.

Note also that, 
if $\d = \d_n \to 0$ sufficiently slowly, the upper bound can be made 
to grow more slowly that $\hPhi^2(n)\log n$.
For example, with $c_2(m) = \log\log m$ and therefore $\hPhi(n) \asymp
\log\log n$, one could take $\d_m = 1/\log\log m$, giving an upper bound 
of order 
\[
  O(\log n \log\log n \log\log\log n) = o(\log n\{\log\log n\}^2)\,. 
\]
In general, taking $\d = \d_n$
to be the solution of the equation $\log(1/\d) = \sqrt{c_2^*(n^\d)}$
gives both $\d_n\to0$ and $\var A_n((1-\d_n)\log n) = o(\hPhi^2(n)\log n)$.
Thus, {\it almost\/} up to the time~$\log n$, the 
compensator~$A_n$ behaves very
much like the simpler integral process~$A_n^*$, but the common scale of their 
fluctuations is of smaller order than that 
of~$A_n(\infty)$, which, by Corollary~\ref{ADB-Psi(n)-bds}, has variance
of order $\Psi_2(n) \asymp \hPhi^2(n)\log n$. 

We now turn to approximating~$A_n(\infty)$.  As before, it is enough to consider
$A_n(\t_n)$, which we can write in the form
\eq\label{time-rev}
  A_n(\t_n) \Eq \int_0^{\t_n} \hPhi(ne^{-S_t})\dt 
	  \Eq \int_0^{\t_n} \hPhi(ne^{-S_{\t_n-v}})\dv.
\en
We now define the process~$\hZ_n$ by the equation 
\eq\label{Zhat-def}
\hZ_n(v) \ :=
\begin{cases}
 S_{\t_n-} - v - S_{(\t_n-v)-}\,\ ~{\rm for~}v < \t_n\,,\\
	 S_{\t_n-} - \t_n\, ~~~~~~~~~~~~~~~~{\rm for~} v \ge \t_n\,,
\end{cases}
\en
and we look for a suitable
approximation to~$A_n(\t_n)$ when the paths of~$\hZ_n$ are `nice'.

Very much as before, we define good events, for $\varphi,\psi>0$,
\eqa
  \hB_1(\f,n) &:=& 
	  \Blb \sup_{0 \le v \le 2\log n} (v\vee1)^{-1/2}|\hZ_n(v)| \le \f \Brb;
		\label{Bhat1-def} \\
	\hB_2(\ps,n) &:=& \{\log n - S_{\t_n-} \le \ps\};\label{Bhat2-def} \\
	\hB_3(T,\ps,n) &:=& \Blb\int_0^T v^{-1}|\hZ_n(v)|\dv \le \ps\Brb,
	  \label{Bhat3-def}
\ena	
whose probabilities we wish to show are large.  The next two lemmas
make this precise;
we recall the definition~\Ref{B0-def} of the event~$B_0(n)$.

\begin{lemma}\label{Bhat}
For any $T,\f,\ps > 0$, we have
\eqs
  \pr[\hB_1^c(\f,n) \cap B_0(n)] &\le& \s^2\f^{-2}(72 + 29\log\log n);\\
	\pr[\hB_3^c(T,\ps,n) \cap B_0(n)] &\le& 6\s\ps^{-1}\sqrt T.
\ens
\end{lemma}

\begin{proof}
In order to make the calculations, it is convenient to exploit the
explicit It{\^o} construction of the process~$S$ \cite[Proposition 1.3]{BertoinSub}.  For~$H$ a Poisson point
process on $\re_+^2$ with intensity measure ${\rm d}t\,\nu_0({\rm d}x)$   
we can define
\[
  S_t \ :=\ \int_{]0,t]\times\re_+} x\,H({\rm d}t\dx);\qquad
    S_t\ut \ :=\ \int_{]0,t]\times\re_+} x^2\,H({\rm d}t\dx),
\]
$S$ being a copy of our original subordinator.
We also define
the family of random point measures~$\mu_t$ on~$\re_+$ by

$$\mu_t[0,x] \ := H(]0,t[\,\times [0,x]).$$
We then define the 
family of $\s$-fields
\[
  \ff_{-t} \ :=\ \s\{\m_t,\,H|_{[t,\infty]\times\re_+}\},\quad t\ge0,
\]
so that $\ff_s \subset \ff_{s'}$ whenever $s \le s' \le 0$.
Then direct calculations show that the processes $(M\ul(t),\,t>0)$, $l=1,2,3$, are
reversed martingales  with respect to the filtration $\{\ff_s,\, s < 0\}$,
with means $1$, $\s^2$ and zero, respectively, where
\[
  M\ui(t) := t^{-1}S_{t-}, \quad M\ut(t) := t^{-1}S_{t-}\ut  \quad\mbox{and}\quad
    M\uh(t) \ :=\ (t^{-1}S_{t-} - 1)^2 - t^{-2}S_{t-}\ut .
\]
Thus it is immediate from the optional 
sampling theorem that
\eq\label{S-means}
  \ex\{\t_n^{-1}S_{\t_n-}\} \Le 1;\qquad
	  \ex\{\t_n^{-1}S_{\t_n-}\ut\} \Le \s^2.
\en
It also follows that $\ex M\uh(t\vee\t_n) = 0$ for any $t > 0$,
which, taking $t = \half\log n$, implies that
\eq\label{S-squared}
  \ex\Blb (\t_n^{-1}S_{\t_n-} - 1)^2 \bone\{\t_n \ge \half\log n\} \Brb
	  \le 2\s^2/\log n.
\en			 		
Furthermore, for $v < \t_n$, the equality $\ex\{M(\t_n-v) \giv \ff_{-\t_n}\}
= W(\t_n)$ a.s.\ also implies that, for such~$v$,
\eq\label{MG-square}
  \ex\{U_n(v)^2 \giv \ff_{-\t_n}\} 
	   \Eq v\,\frac{\t_n^{-1}S_{\t_n-}\ut}{\t_n(\t_n-v)}\,,
\en
where 
\[
  U_n(v) \ :=\ (\t_n-v)^{-1}S_{(\t_n-v)-} - \t_n^{-1}S_{\t_n-}\,.
\]

We thus have the expression
\eq\label{Zhat-def2}
  \hZ_n(v) \Eq (v\wedge\t_n)\{\t_n^{-1}S_{\t_n-} - 1\}  
	  - (\t_n-v)_+U_n(v),\quad v\ge0,
\en
as an alternative representation for~$\hZ_n$, in addition to~\Ref{Zhat-def}.
Taking expectations conditional of~$\ff_{-\t_n}$, we thus obtain
\eqs
  \ex\{|\hZ_n(v)| \giv \ff_{-\t_n}\}	&\le& (v\wedge\t_n)|\t_n^{-1}S_{\t_n-} - 1|
	  + (\t_n-v)_+\ex\{|U_n(v)| \giv \ff_{-\t_n}\}  \\
	&\le& v|\t_n^{-1}S_{\t_n-} - 1| + \sqrt{v}\,\sqrt{\t_n^{-1}S_{\t_n-}\ut},
\ens
the last inequality from~\Ref{MG-square}. Multiplying by 
$\bone\{\t_n \ge \half\log n\}$ and taking expectations thus yields
\eq\label{Ex-bnd}
	\ex(|\hZ_n(v)|\bone\{\t_n \ge \half\log n\}) 
	\Le v\s\sqrt{2/\log n} + \s\sqrt v \Le 3\s\sqrt v,
\en
for $0 \le v \le 2\log n$, in view of \Ref{S-means} and~\Ref{S-squared}.
The second inequality now follows from Markov's inequality, because
\[
  \ex\Blb \int_0^T v^{-1}|\hZ_n(v)|\dv \bone\{B_0(n)\}\Brb \Le 3\s\int_0^T v^{-1/2}\,{\rm d}v.
\]

It also follows from~\Ref{Zhat-def2} that, for any~$\f > 0$ and for $v < \t_n$,	
\[
  \{|\hZ_n(v)| > \f\sqrt{v\vee1}\} \subset 
	  \{ |\t_n^{-1}S_{\t_n-} - 1| > \half\f v^{-1/2}\}
	  \cup \{(\t_n-v)|U_n(v)| > \half\f v^{1/2}\}\,.
\]
The first event happens for some $v < \t_n$ only if	
$|\t_n^{-1}S_{\t_n-} - 1| > \half\f\t_n^{-1/2}$, and the probability of this
happening on the event $B_0(n)$ is at most
\eqs
  \lefteqn{\pr\Bigl[|\t_n^{-1}S_{\t_n-} - 1|\bone\{\t_n \ge \half\log n\} 
	   > \half\f(2\log n)^{-1/2}\Bigr]}\\
	&&\le\ \frac{2\s^2}{\log n}\cdot\frac{2\log n}{\f^2} \Eq 4\s^2\f^{-2},
	\hspace{2cm}
\ens
by~\Ref{S-squared}.	For the second, using Kolmogorov's inequality much 
as in the proof of Lemma~\ref{B-probs}, for~$r\ge1$ such that 
$2^r \le \half\t_n$, we have
\eqs
  \lefteqn{\pr\Bigl[\sup_{2^{r-1} \le v \le 2^r} v^{-1/2}(\t_n-v)|U_n(v)| 
	   > \half\f \giv \ff_{-\t_n}\Bigr]} \\
	&&\le\ \pr\Bigl[\sup_{0 \le v \le 2^r} |U_n(v)| 
	   > \half\t_n^{-1}\f 2^{(r-1)/2}\giv \ff_{-\t_n}\Bigr] \\
	&&\le\ \frac{4\ex\{|U_n(2^r)|^2 \giv \ff_{-\t_n}\}\t_n^2}{\f^2\,2^{r-1}} \\
  &&\le\  16\f^{-2}\t_n^{-1}S_{\t_n-}\ut,
\ens
using~\Ref{MG-square}.  Adding over all such~$r$, and including the
$v$-intervals $]0,1[$ and $]2^r,\half\t_n]$, it follows that
\eq\label{Zhat-bnd1}
  \pr\Bigl[\sup_{0 \le v \le \t_n/2} v^{-1/2}(\t_n-v)|U_n(v)| 
	   \bone\{B_0(n)\} > \half\f\Bigr\}	
		 \Le 20\s^2\f^{-2}(1 + \lceil \log_2\log n \rceil),
\en
from~\Ref{S-means}.  For $\half\log n < v < \t_n$, we use~\Ref{Zhat-def}
to give
\[
   \hZ_n(v) \Eq  Z_{\t_n-} - Z_{(\t_n-v)-},
\]
so that
\eqs
	&&B_0(n) \cap 
	  \Blb \sup_{\shalf\t_n \le v\le \t_n} v^{-1/2}|\hZ_n(v)| > \f \Brb
	\ \subset\ \Blb \sup_{0 \le u \le 2\log n}|Z_u| > \quarter\f\sqln \Brb,
   \ens
the latter event, by Kolmogorov's inequality, having probability at most
$32\s^2\f^{-2}$.  Finally, again by Kolmogorov's inequality,
\[
  \pr[B_0(n)^c] 
	  \Le \pr\Bigl[\sup_{0 \le u \le 2\log n}|Z_u| > \half\log n\Bigr]
	  \Le 8\s^2/\log n.
\]
From these last two bounds and from~\Ref{Zhat-bnd1}, the lemma follows.
\hfill$\Box$		
\end{proof}

\msk
\begin{lemma}\label{last-jump}
If $\ps_n \to \infty$, then $\lnti\pr[\hB_2^c(\ps_n,n)] = 0$. 
\end{lemma}

\begin{proof}
Simply note that $\hB_2^c(x,n) \subset \hB_2^c(y,n)$ whenever $x > y$, 
so that then
\[
  \pi_n(x) \ :=\ \pr[\hB_2^c(x,n)] \Le \pi_n(y),
\]
and that
\[
  \lnti\pi_n(x) \Eq 
\frac{\int_x^\infty N_0(u)\du} 
    {\int_0^\infty N_0(u)\du} \ =:\ \pi(x),
\]
by the  renewal theorem \cite[p. 99]{Bertoin}, 
with $\lim_{x\to\infty}\pi(x) = 0$.
Hence, given $\e > 0$, pick~$x$ so that $\pi(x) < \e/2$, and then~$n_x$
such that $\pi_n(x) \le \e$ and $\ps_n \ge x$ for all $n \ge n_x$;
it then follows that $\pi_n(\ps_n) \le \pi_n(x) \le \e$ for all $n \ge n_x$.
\hfill$\Box$
\end{proof}

\bigskip
With these preparations, we are now in a position to approximate the
behaviour of~$A_n(\t_n)$, and indeed of the whole process
$A_n(t\wedge \t_n)$.

\begin{theorem}\label{VSGPhi}
Suppose that Assumptions {\rm A1--A2} and~\Refrm{L-cond6} hold.
Fix $\a=1/8$, $\b=1/4$, and set $v_n := 4\log^{2\a}n$.  Then,
on the event 
\[
  B_0(n)\cap\hB_1(\log^\a n,n)\cap\hB_2(\log^\b n,n)
	  \cap\hB_3(2\log n,\sqrt{c_2^*(e^{v_n})\log n},n),
\]
it follows that
\[
  \Bl\hPhi(n)\sqln\Br^{-1}\sup_{t\ge0}
	   \Blm A_n(t\wedge\t_n) - \int_{(\t_n-t)_+}^{\t_n} \hPhi(e^v)\dv\Brm	
		 \ \to\ 0.
\]
Furthermore,
\[
  \lnti\pr[B_0(n)\cap\hB_1(\log^\a n,n)\cap\hB_2(\log^\b n,n)
	  \cap\hB_3(2\log n,\sqrt{c_2^*(e^{v_n})\log n},n)] = 1.
\]		
Here, $c_2^*(m)$ is defined as in {\rm Lemma \ref{T1-new}}.

\end{theorem}		 	 

\begin{proof}
Recalling~\Ref{time-rev}, we can write
\eqa
  A_n(t\wedge\t_n) &=& \int_0^{(t\wedge\t_n)} \hPhi(ne^{-S_u})\du\non\\
  &=&\ \int_{(\t_n-t)_+}^{\t_n} \hPhi(e^v)\dv
	  + \int_{(\t_n-t)_+}^{\t_n} \Blb \hPhi(e^{v+D_n}) - \hPhi(e^v)\Brb\dv\non\\
	&&\quad\mbox{}
	  + \int_{(\t_n-t)_+}^{\t_n}
		\Blb \hPhi(e^{v+D_n+\hZ_n(v)}) - \hPhi(e^{v+D_n}) \Brb\dv,
	\label{Vslow-1}
\ena
where $D_n := \log n - S_{\t_n-} \ge 0$.  The second of the
integrals in~\Ref{Vslow-1} is nonnegative, and no larger than
\eqa
  \int_{0}^{\t_n} \Blb \hPhi(e^{v+D_n}) - \hPhi(e^v)\Brb\dv
	  &\le& \int_{\t_n}^{\t_n+D_n}\hPhi(e^v)\dv\non\\
	&\le& \hPhi(n^3)\log^\b n,\label{2nd-term}
\ena
on the event $B_0(n)\cap \hB_2(\log^\b n,n)$.  Note also that, for
any $r\ge1$ and~$n$ such that $c_2^*(n) \ge 1$,
\eqa
  1 &\le& \frac{\hPhi(n^r)}{\hPhi(n)} 
	  \Eq \exp\Blb \int_n^{n^r} \frac{{\rm d}y}{yL(y)} \Brb\non\\
	&\le& \exp\Blb (\log\log(n^r) - \log\log n)/c_2^*(n) \Brb
	  \Eq r^{1/c_2^*(n)} \Le r; \label{Phi-ratio}
\ena
hence, from~\Ref{2nd-term}, the second of the integrals in~\Ref{Vslow-1} 
is of smaller
order than $\hPhi(n)\sqln$ on the event $B_0(n)\cap \hB_2(\log^\b n,n)$.

To control the third of the integrals in~\Ref{Vslow-1}, we bound
\eq\label{3rd-term}
  \int_{0}^{\t_n}
		\Blm \hPhi(e^{v+D_n+\hZ_n(v)}) - \hPhi(e^{v+D_n}) \Brm\dv\,.
\en		
On $\hB_1(\log^\a n,n)$, we have $v^{-1}|\hZ_n(v)| \le 1/2$ if $v\ge v_n$.
So split the range of the integral into $0 < v \le v_n$ and $v_n \le v
\le \t_n$.  In the lower range, on $\hB_1(\log^\a n,n)\cap \hB_2(\log^\b n,n)$,
the exponents $v + D_n$ and $v+D_n + \hZ_n(v)$ are bounded above by
\[
  v_n + \log^\b n + \log^\a n\sqrt{v_n} \Le 7\log n,
\]
implying, together with~\Ref{Phi-ratio}, that~\Ref{3rd-term} is
bounded above by $7\hPhi(n)v_n$ for all~$n$ large enough, and this
is $o(\hPhi(n)\sqln)$ by choice of~$\a$.  In the upper range,	
we can apply Lemma~\ref{keylemma} to $\hPhi(e^{v+D_n+\hZ_n(v)})$, very 
much as in the proof of Theorem~\ref{T1}, because here 
$|\hZ_n(v)| \le \half(v+D_n)$. The quantity $\hX(n,v)$,
analogous to~$X(n,t)$ of \Ref{Xnt-def}, is bounded for $v\ge v_n$ by
\[
   \hX(n,v) \Le \frac{\k|\hZ_n(v)|^2}{c_2(e^{v+D_n})v(v+D_n)}
	   \Le \frac{\k\log^{2\a}n}{vc_2^*(e^{v_n})} \Le \frac\k{c_2^*(e^{v_n})},
\]
and
\[
  \frac{|\hZ_n(v)|}{L(e^{v+D_n})} \Le \frac{\log^\a n}{c_2^*(e^{v_n})\sqrt v}
	  \Le \frac1{2c_2^*(e^{v_n})},
\]				 
giving
\eqa
   \lefteqn{\int_{v_n}^{\t_n}\Blm \hPhi(e^{v+D_n+\hZ_n(v)}) 
	    - \hPhi(e^{v+D_n}) \Brm\dv}\non\\
	 &\le& 
	   \int_{v_n}^{\t_n} \hPhi(e^{v+D_n})\frac{|\hZ_n(v)|}{vc_2^*(e^v)}\dv
	    +	\lambda\Bl \frac{\hPhi(e^{\t_n+D_n}) \log^{2\a}\!n}
			    {c_2^*(e^{v_n})}\int_{v_n}^{\t_n} v^{-1}\dv	\Br,\phantom{HHH}
					\label{Vslow-2}
\ena
for some positive constant $\lambda < \infty$.
On $B_0(n)\cap \hB_2(\log^\b n,n)$, and from~\Ref{Phi-ratio}, 
we have $\hPhi(e^{\t_n+D_n}) \le 3\hPhi(n)$ for all~$n$ large enough,
so that the second term in~\Ref{Vslow-2} is of order $o(\hPhi(n)\sqln)$.
The first term is bounded 
on
\[
  B_0(n)\cap\hB_2(\log^\b n,n)\cap\hB_3(2\log n,\sqrt{c_2^*(e^{v_n})\log n},n) 
\]
by 
\[
  \hPhi(n^3)\sqln\,/\sqrt{c_2^*(e^{v_n})} \Eq o(\hPhi(n)\sqln),
\]
again by~\Ref{Phi-ratio}.  This completes the proof of the main
statement.  The final assertion follows from Lemmas \ref{B-probs},
\ref{Bhat} and~\ref{last-jump}.
\hfill$\Box$
\end{proof}

\section{Approximation theorems}
We can now build on the results of the previous sections to derive
central limit approximations for~$\kk_n$.  The starting point is
the functional central limit theorem for the L\'evy process itself.
Defining the process~$W_m$ by $W_m(t) := \s^{-1} m^{-1/2}Z_{mt}$,
it follows that
\eq\label{Z-CLT}
  W_m\ \to_d\ W \quad\mbox{in}\ D_1[0,\infty[\quad\mbox{as}\ m \to \infty,
\en
where~$W$ is standard Brownian motion and $D_1[0,\infty[\,$ denotes 
the space of c\`adl\`ag functions	
$x\colon [0,\infty[\, \to \re$ satisfying $\lim_{t\to\infty}t^{-1}x(t)
= 0$, endowed with the metric $\r_1(x,y) := \sup_{t\ge0}(t\vee1)^{-1}
|x(t)-y(t)|$ (M{\"u}ller~\cite{Mueller}, Satz~1). 
As a consequence of the central limit theorem for the renewal 
processes \cite[Section XI.5]{Feller},
it also follows that
\eq\label{t_n-CLT}
  U_n\ :=\ (\t_n - \log n)/\{\s\sqln\}\ \to_d\ \nn(0,1) 
	  \quad\mbox{as}\ n \to \infty.
\en

\par
We shall also be interested in approximations which are not given in the form of
limit theorems, but are instead expressed in terms of bounds on a distance
between the distributions of the processes considered, taken here to be the
appropriate bounded Wasserstein distances.  For probability measures $Q$ 
and~$Q'$ on a metric space $(\xx,\r)$, the bounded Wasserstein distance 
$d_{BW}(Q,Q')$ is
defined to be $\sup_{f\in\ww}|\int f\,\rmd Q - \int f\,\rmd Q'|$, where~$\ww$
denotes the bounded Lipschitz functions on~$\xx$:
\[
  \ww\ :=\ \ww_{\xx,\r}\ :=\ \{f\colon \xx\to\re: \|f\| \le 1,\,
	  L(f) \le 1\},
\]
and $L(f) := \sup_{x\neq x'\in\xx} |f(x)-f(x')|/\r(x,x')$.  The distance~$d_{BW}$
metrises
weak convergence in  $(\xx,\r)$ (Dudley \cite{Dudley}, Theorem~8.3).  Note 
also that if, for each $n\ge1$, the random elements $X_n$ and~$Y_n$ of $(\xx,\r)$ are
on the same probability space, then $d_{BW}(\law(X_n),\law(Y_n)) \to 0$ if,
for each $\e>0$, $\pr[\r(X_n,Y_n) > \e] \to 0$.  If~$\xx$ is the space
$D_1[0,\infty[\,$ defined above, we shall refer to $\ww^1$ and~$d_{BW}^1$; 
if~$\xx$ is the space $D_0[0,\infty[\,$ of c\`adl\`ag functions	
$x\colon [0,\infty[\, \to \re$ having finite limits as $t\to\infty$, endowed 
with the metric $\r_0(x,y) := \sup_{t\ge0}|x(t)-y(t)|$, we shall refer to $\ww^0$ 
and~$d_{BW}^0$, and, if $\xx=\re$, we shall write~$d_{BW}^\re$.

\subsection{Moderate growth}
We begin with a setting of moderate growth, in which $L(n)
\asymp \log n$, so that~\Ref{L-cond1} is in force. 
In order to
describe the behaviour of~$\kk_n$, we first define a centred
and normalized version~$\kk_n\ui$ of the process by
\[
  \kk_n\ui(u) \ :=\ \Bl\hPhi(n)\sqln\Br^{-1}\Blb
	  \kk_n(u\log n) - \log n\int_0^{(u\wedge1)}\hPhi(n^{1-v})\dv\Brb,
\]
whose distribution we approximate by that of~$Y_n\ui$, where
\[
  Y_n\ui(u) \ :=\ \s\int_0^{(u\wedge1)} h_n\ui(v)\, W(v)\dv,
\]
with
\[
  h_n\ui(u)\ :=\ \frac{\hPhi(n^{1-u})\log n}{\hPhi(n)L(n^{1-u})}.
\]
Note that $h_n\ui(u) \ge 0$ for all~$u$, and that, from~\eqref{gnt}, 
\[
  \int_0^1 h_n\ui(u) \du \Eq 1 - \hPhi(1)/\hPhi(n) \Le 1\,~~~~~~~(n\geq 1).
\]

\begin{theorem}\label{MG-T1}
  If Assumptions {\rm A1--A2} hold, and $L(n) \asymp \log n$, 
then
\[
  d_{BW}(\law(\kk_n\ui),\law(Y_n\ui))\ \to\ 0 \quad\mbox{as }
	  \nti\,.
\]

\end{theorem}

\begin{proof}
We begin by writing
\eqs
  \lefteqn{\kk_n\ui(u) - \tY_n\ui(u)
	\Eq \ r_n\Bigl\{
	  (\kk_n(u\log n) - A_n(u\log n))}\\ 
  &&\mbox{}\ + (A_n(u\log n) - A_n(\{u\log n\}\wedge\t_n))
		+ (A_n(\{u\log n\}\wedge\t_n) - A_n^*(u\log n))\Bigr\}\,,
\ens
where $r_n := \Bl\hPhi(n)\sqln\Br^{-1}$ and
\eqs
  \tY_n\ui &:=& A_n^*(u\log n) - \log n \int_0^{u\wedge1} \hPhi(n^{1-v})\dv\\
	&=&  -\s\int_0^{(u\wedge1)} h_n\ui(v)\, W_{\log n}(v)\dv.
\ens
Now we have $r_n\sup_{u\ge0}|\kk_n(u\log n) - A_n(u\log n)| \to_p 0$ by Lemma~\ref{KAlink}
and Corollary~\ref{ADB-Psi(n)-bds}, then
$r_n\sup_{u\ge0}|A_n(u\log n) - A_n(\{u\log n\}\wedge\t_n)| \to_p 0$ by 
Lemma~\ref{tau-truncation}, and finally, by Theorem~\ref{T1}, 
$r_n\sup_{u\ge0}|A_n(\{u\log n\}\wedge\t_n) - A_n^*(u\log n)| \to_p 0$. 
Hence it follows that
\[
  d_{BW}^0(\law(\kk_n\ui),\law(\tY_n\ui))\ \to\ 0 \quad\mbox{as }
	  \nti.
\]

To conclude the proof, we now need to show that 
$\sup_{f\in \ww^0}|\ex f(\tY_n\ui) - \ex f(Y_n\ui)|\to0$ as~$\nti$.  To do so,
for any $f\in\ww^0$, 
define $f_n\colon D_1[0,\infty[\,\to\re$ by $f_n(w) := f(H_n(w))$,
where $H_n(w)(u) := \int_0^{(u\wedge1)} h_n\ui(v)w(v)\dv$.  Note that, for 
$w,w'\in D_1[0,\infty[\,$ and any $u\ge0$,
\eqs
  |H_n(w)(u) - H_n(w')(u)| &=& \Blm \int_0^{(u\wedge1)} h_n\ui(v)(w(v) - w'(v))\dv \Brm\\
  &\le& \r_1(w,w')\int_0^1 h_n\ui(v)\dv \Le \r_1(w,w').
\ens
Hence, for any $f\in\ww^0$, it follows that $f_n \in \ww^1$, and hence that
\[
  |\ex f(\tY_n\ui) - \ex f(Y_n\ui)| \Eq |\ex f_n(W_{\log n}) - \ex f_n(W)| 
  \Le d_{BW}^1(\law(W_{\log n}),\law(W)). 
\]   
The theorem now follows from~\Ref{Z-CLT}.
\hfill$\Box$
\end{proof}

\begin{theorem}\label{MG-T2}
Under the assumptions of {\rm Theorem~\ref{MG-T1}}, if in addition
$L(n) \sim \g\log n$ 
for some $0 < \g < \infty$, then
\[
  \kk_n\ui\ \to_d\ Y\ui\quad\mbox{in } D_0[0,\infty[\,\quad\mbox{as }
	  \nti,
\]
where
\[
  Y\ui(u) \ :=\ \s\int_0^{(u\wedge1)} \g^{-1} (1-v)^{(1-\g)/\g} W(v)\dv.
\]		
\end{theorem}

\begin{proof}
If $L(n) \sim \g\log n$, then $h_n\ui(u) \to \g^{-1} (1-u)^{(1-\g)/\g}$
uniformly in $0 \le u \le 1-\d$, for any $\d > 0$;
furthermore,
\[
  \limsup_{\nti}\int_{1-\d}^1 h_n\ui(v)\dv 
    \Eq \limsup_{\nti}\{\hPhi(n^\d)-\hPhi(1)\}/\hPhi(n) \Le \d^{\g'}
\]
for any $\g' < \g$, and $\int_{1-\d}^1 h\ui(v)\dv = \d^\g$. Hence
\[
  \ex\Blb \sup_{u\ge0}|Y_n\ui(u) - Y\ui(u)| \Brb \le \sqrt{\frac2{\pi}}
    \int_0^1 |h_n\ui(v) - h(v)|\dv \ \to\ 0,
\]
proving the theorem.        
\hfill$\Box$
\end{proof}

\vskip0.5cm
\noindent
{\bf Examples.}
Suppose, for some $0 < \g < \infty$, that~$\tS$ is a subordinator such that
$\tL(n) \sim \g\log n$ and $\tPhi(n)\sim c\log^{1/\gamma} n$; as at the end
of Section~\ref{variance}, we do not assume that $\mm = \ex\tS_1$ takes the
value~$1$, and we write $\t^2 = \var\tS_1$.
Theorem \ref{MG-T2} entails a gaussian limit for $(K_n - \mu_n)/\s_n$, with 
$$
  \mu_n\sim \log n\int_0^1\hPhi(n^{1-v})\dv\ \sim\ 
    \frac{c\log^{1+1/\gamma} n}{\mm(1+1/\gamma)}
$$
and
\[    
    \s_n^2\ \sim\ \hPhi^2(n)\log n\,\s^2\, \var \left\{
      \int_0^1 \gamma^{-1}(1-v)^{1/\gamma-1}W(v){\rm d}v\right\}
      \ \sim\ \frac{c^2\t^2 \log^{1+2/\gamma} n}{\mm^3(1+2/\gamma)}\,, 
\]
where, as before, $\hPhi(n) = \mm^{-1}\tPhi(n)$ and $\s^2 = \mm^{-1}\t^2$.
Note also that $\mu_n \sim \Psi(n) \sim \mm^{-1}\tPsi(n)$ and that 
$\s_n^2 \sim \s^2\Psi_2(n)\sim \t^2\mm^{-3}\tPsi_2(n)$, as is to be expected.

\par For the classical gamma subordinator \cite[p. 73]{Bertoin}, scaling so that
$\ex S_1 = 1$, we have $\nu_0({\rm d}x)=\th e^{-\theta x}{\rm d}x/x$,
$\Phi_0(n)=\th\log (1+n/\theta)$, and $\s^2 = 1/\th$.
Hence the CLT in \cite{GPYII} agrees with Theorem \ref{MG-T2}. Note that 
one parameter $\theta>0$ is enough, since, for the L{\'e}vy measure $a\nu_0$, the 
distribution of $K_n$ does not depend on the scale parameter~$a$.

\par In the case $\gamma=1$, Theorem \ref{MG-T2}  covers a somewhat larger family
of gamma-like subordinators than that  considered in \cite{GPYII}. The extension is that 
the condition of exponential decay for $N_0(x)$ as $x\to\infty$ required in \cite{GPYII} 
is replaced now by a weaker condition $\sigma^2<\infty$. The constraints on the 
behaviour of $N(x)$ at $x\to 0$ are also slightly weaker here.

\subsection{Fast growth}
We now turn to the setting in which $L(n) \to \infty$ but $L(n)/\log n \to 0$; hence
$\Phi$ grows faster than any power of the logarithm. In order to apply the previous
theorems, we need to suppose either that~\Ref{L-cond1} is in force,
albeit with $L(n) = o(\log n)$, or that~\Ref{L-cond4} holds. The
analogue of $\kk_n\ui$ is now~$\kk_n\ut$, defined by
\[
  \kk_n\ut(u) \ :=\ \Bl\hPhi(n)\sqLn\Br^{-1}\Blb
	  \kk_n(uL(n)) - L(n)\int_0^{(u\wedge l_n)}\hPhi(ne^{-vL(n)})\dv\Brb,
\]
where $l_n := \log n/L(n)$.  Here,
we approximate the distribution of~$\kk_n\ut$ by that of~$Y\ut$, where
\[
  Y\ut(u) \ :=\ \s\int_0^{u} e^{-v}\, W(v)\dv.
\]

\begin{theorem}\label{FG-T1}
If Assumptions {\rm A1--A2} hold, and $L(n)/ \log n \to 0$, with
either \Refrm{L-cond1} or~\Refrm{L-cond4} satisfied, then
\[
  d_{BW}(\law(\kk_n\ut),\law(Y\ut))\ \to\ 0 \quad\mbox{as }
	  \nti,
\]
where $d_{BW}$ is as before.
\end{theorem}

\begin{proof}
If \Ref{L-cond1} is satisfied, we argue as in the proof of Theorem~\ref{MG-T1} 
to show that
\eq\label{ADB:limit-2.1}
  \sup_{u\ge0} |\kk_n\ut(u) - \tY_n\ut(u)|\ \to_p\ 0 \quad\mbox{as }
	  \nti,
\en
where 
\[
  \tY_n\ut(u) \ :=\ \s\int_0^{(u\wedge l_n)} h_n\ut(v) W_{L(n)}(v)\dv,
\]
and
\[
  h_n\ut(u)\ :=\ \frac{\hPhi(ne^{-uL(n)})L(n)}{\hPhi(n)L(ne^{-uL(n)})}.
\]
If \Ref{L-cond4} is satisfied, \eqref{ADB:limit-2.1} is still true, using 
Theorem~\ref{VFGPhi} in place of Theorem~\ref{T1} in the proof.
Once again, $h_n\ut(u) \ge 0$ for all~$u$, and  
\[
  \int_0^\infty h_n\ut(u) \du \Eq 1.
\]

The next step is to approximate $\tY_n\ut$ by~$Y_n\ut$, where
\[
   Y_n\ut(u) \ :=\ \s\int_0^{u} e^{-v} W_{L(n)}(v)\dv.
\] 
Here, it is immediate that
\eq\label{ADB:limit-2.3}
  \ex\Blb \sup_{u\ge0}\s^{-1}|\tY_n\ut(u) - Y_n\ut(u)| \Brb
    \Le \int_0^{l'_n} |h_n\ut(v) - e^{-v}| \sqrt{v}\dv
    + \int_{l'_n}^\infty \{h_n\ut(v) + e^{-v}\} \sqrt{v}\dv
\en
for any $l'_n \le l_n$; we take $l'_n = \min\{l_n^{1/2},\tfrac1{2(\k\vee1)}l_n\}$.
Now, from~\eqref{bounds}, for $0 \le v\le l'_n$, we have
\[
  \Blm \frac{L(n)}{L(ne^{-uL(n)})} - 1 \Brm \Le \l_1 l_n^{-1/2},
\]
for some $\l_1 < \infty$, and hence, by Corollary~\ref{ADB-cross}, that 
\[
   \sqrt{v}\Blm h_n\ut(v) - \Phi^*(n,v)\Brm
     \Le \l_2\sqrt{v}e^{-v/2} l_n^{-1/2},
\]
where 
\[
  \Phi^*(n,v) \ :=\ \frac{\hPhi(ne^{-vL(n)})}{\hPhi(n)}.
\]  
Then Lemma~\ref{keylemma} gives
\[
  \sqrt v \Blm  \Phi^*(n,v) - e^{-v}\Brm 
    \Le \l_3\{\exp(\k v^2/l_n) - 1\}\sqrt v e^{-v} \Le \l_4 l_n^{-1}v^{5/2}e^{-v},
\]
for different constants $\l_3,\l_4$, again in $0 \le v\le l'_n$.  Hence it follows
that
\eq\label{ADB:limit-2.2}
  \lnti \int_0^{l'_n} |h_n\ut(v) - e^{-v}| \sqrt{v}\dv \Eq 0.
\en
 
It is also immediate that $\lnti \int_{l'_n}^\infty e^{-v} \sqrt{v}\dv  = 0$.
Hence, to show that the right hand side of~\eqref{ADB:limit-2.3} is small in the
limit, it remains only to consider
\eq\label{ADB:limit-2.4}
  \int_{l'_n}^\infty h_n\ut(v) \sqrt{v}\dv 
    \Eq \sqrt{l'_n}\Phi^*(n,l'_n)
    + \int_{l'_n}^\infty \half v^{-1/2}\Phi^*(n,v)\dv.
\en
Here, the first term tends to zero as $\nti$ by Corollary~\ref{ADB-cross}, as does
\[
  \int_{l'_n}^{\a_n} v^{-1/2}\Phi^*(n,v)\dv \Le \int_{l'_n}^{\a_n} v^{-1/2}e^{-v/2}\dv,
\]      
where $\a_n := l_n/\{2(\k\vee1)\}$.  Then, splitting the remaining integral at~$2l_n$,
we have
\[
  \int_{\a_n}^\infty v^{-1/2}\Phi^*(n,v)\dv \Le 2l_n \exp\{-\a_n/2\} + 
    \frac1{n\hPhi(n)L(n)}\,,
\]
the final term following from~\eqref{Phi-bnd-1}. Combining these bounds, we have 
now also shown
that $\lnti \int_{l'_n}^\infty h_n\ut(v) \sqrt{v}\dv = 0$; hence, from
\eqref{ADB:limit-2.1} and~\eqref{ADB:limit-2.3}, it follows that
\[
   d_{BW}^0(\law(\kk_n\ut),\law(Y_n\ut)) \ \to\ 0 \quad\mbox{as}\ \nti.\   
\]

Finally, if $H\colon D_1[0,\infty[\, \to D_0[0,\infty[\,$ is
defined by $H(w)(u) := \int_0^u e^{-v}w(v)\dv$ and~$f$ is in~$\ww^0$, 
then $(1+e^{-1})^{-1}f\circ H \in \ww^1$,
from which $d_{BW}^0(\law(Y_n\ut),\law(Y\ut)) \to 0$ as
$\nti$ follows immediately, and the theorem is proved.
\endpf
\end{proof}

\subsection{Slow growth}
If~$\hPhi$ grows very slowly to infinity, with $L(n)/\log n
\to\infty$, the arguments culminating in Theorem~\ref{VSGPhi}
show that the key quantity describing the process~$\kk_n$ is the
family of integrals
\[											
  \int_{(\t_n-t)_+}^{\t_n} \hPhi(e^v)\dv,\quad t\ge0.
\]
Here, the randomness enters only through the hitting time~$\t_n$,
which is asymptotically normally distributed, as recorded in~\Ref{t_n-CLT}.
The process thus has a quite different qualitative behaviour to
that of the previous cases.

Since $\t_n$ takes values fairly close to~$\log n$, it makes sense
to describe the random behaviour of~$\kk_n(t)$ by first	
subtracting $\int_{(\log n-t)_+}^{\log n} \hPhi(e^v)\dv$, and
then dividing by $\hPhi(n)\sqln$.  This leads us to define the
process~$\kk_n\uh$ for $t\ge0$ by
\[
  \kk_n\uh(t) \ :=\ \Bl \hPhi(n)\sqln \Br^{-1}
	  \Blb\kk_n(t) - \int_{(\log n-t)_+}^{\log n} \hPhi(e^v)\dv\Brb.
\]
Then, defining $G_n\colon \re\to D_0[0,\infty[\,$ by
\eq\label{Y3-def}
  G_n[u](t) \ :=\ \s u - \Bl \hPhi(n)\sqln \Br^{-1}
	  \int_{(\log n-t)_+}^{(\log n - t + \s u\sqln)_+} \hPhi(e^v)\dv,
	  \quad t\ge0,
\en
for each $u\in\re$,
we define our approximating process to be ~$Y_n\uh = G_n[U]$,
where~$U$ is a standard normal random variable.

\begin{theorem}
Suppose that Assumptions {\rm A1--A2} and~\Refrm{L-cond6} hold.  Then
it follows that
\[
  d_{BW}^0(\law(\kk_n\uh),\law(Y_n\uh))\ \to\ 0\quad\mbox{as }\nti.
\]
\end{theorem}

\begin{proof}
Once again, we combine Lemmas \ref{tau-truncation} and~\ref{KAlink} and
Corollary~\ref{ADB-Psi(n)-bds}, this time with Theorem~\ref{VSGPhi}, 
showing that
\[
  d_{BW}^0(\law(\kk_n\uh),\law(\tY_n\uh))\ \to\ 0 \quad\mbox{as }
	  \nti,
\]
where 
\eqs
  \tY_n\uh(t) &:=& \Bl \hPhi(n)\sqln \Br^{-1}
	  \Blb \int_{(\t_n-t)_+}^{\t_n} \hPhi(e^v)\dv - 
		  \int_{(\log n-t)_+}^{\log n} \hPhi(e^v)\dv\Brb\\
	&=& \Bl \hPhi(n)\sqln \Br^{-1}
	  \Blb \int_{\log n}^{\t_n}  \hPhi(e^v)\dv - 
	  \int_{(\log n-t)_+}^{(\t_n-t)_+} \hPhi(e^v)\dv\Brb.
\ens
Now, from~\eqref{Phi-ratio}, we have
\[
  \bone\{B_0(n)\}\Bl \hPhi(n)\sqln \Br^{-1}\Blm \int_{\log n}^{\t_n}  \hPhi(e^v)\dv 
    - \int_{\log n}^{\t_n}  \hPhi(n)\dv \Brm \Le |(2^{1/c_2^*(\sqrt n)} - 1)U_n|,
\]
and this tends to zero as $\nti$ because of~\eqref{t_n-CLT} together with $\lnti
c_2^*(n) = \infty$.  On the other hand, $\lnti\pr[B_0(n)]  = 0$, by Lemma~\ref{B-probs}.
Hence $d_{BW}^0(\law(\tY_n\uh),\law(G_n[U_n])) \to 0$ as $\nti$.

To complete the proof, we just have to show that the distributions of $G_n[U_n]$
and $G_n[U]$ are close. To do so, we first define $\tG_n\colon \re\to D_0[0,\infty[\,$ by 
$\tG_n[u] = G_n[u\wedge\s^{-1}\sqrt{\log n}]$.  Then, from the definition of~$G_n$
and from~\eqref{Phi-ratio}, it is immediate that
\eqs
  \sup_{v\ge0}|\tG_n[u'](v) - \tG_n[u](v)| &\le& \s|u'-u| +
    \Bl \hPhi(n)\sqln \Br^{-1} \hPhi(n^2) \s\sqln|u'-u|\\
  &\le& 3\s|u'-u|;
\ens
hence $\lnti\sup_{f\in\ww^0}|\ex f(\tG_n[U_n]) - \ex f(\tG_n[U])| = 0$, in view of
\eqref{t_n-CLT}.  Finally,
\eqs
  \lefteqn{\sup_{f\in\ww^0}\{|\ex f(G_n[U_n]) - \ex f(\tG_n[U_n])| +  |\ex f(G_n[U])
    - \ex f(\tG_n[U])| \}}\\
    &&\Le 2\pr[U_n > \s^{-1}\sqrt{\log n}] + 2\pr[U > \s^{-1}\sqrt{\log n}] \ \to\ 0,
\ens
and the theorem follows.    
\hfill$\Box$
\end{proof}

\bigskip
The process $Y_n\uh$ starts close to zero, and, as indicated by 
Lemmas \ref{T1-new} and~\ref{Var-2}, remains close to zero
until $1- t/\log n$ becomes small.  It reaches its final value
$\s U$ at time $\log n + \s U\sqln$ if $U \ge 0$, and at time
$\log n$ if $U < 0$. 

Its behaviour can also be understood in
terms of the overlapping representation provided under the
condition~\Ref{L-cond2}, when~$c_2(n)$ is allowed to tend to
infinity, but not too fast. Here, the approximation to the
random fluctuations is expressed in terms of the process
\[
  \Bl\hPhi(n)\sqln\Br^{-1}
    \int_0^{T\wedge\log n} \hPhi(ne^{-t})\frac{Z_t}{L(ne^{-t})}\dt,
\]
which at first sight looks very different.  Here, however,
as already observed at the start of Section~\ref{backward},
\[
  \Bl\hPhi(n)\sqln\Br^{-1}
	  \int_0^{(u\log n)\wedge\log n} \hPhi(ne^{-t})\frac{Z_t}{L(ne^{-t})}\dt
\]
is of small order whenever~$u$ is bounded away from~$1$, and even 
for choices of~$u = u(n) \to 1$ such that
$(1 - u(n))c_2(n^{1-u(n)}) \to \infty$.  On the other hand, for~$u$				 
closer to~$1$, the 
remaining contribution is approximately
\eqs
   \lefteqn{\Bl\hPhi(n)\sqln\Br^{-1}Z_{\log n}\int_{u(n)\log n}^{u\log n} 
	    \frac{\hPhi(ne^{-t})}{L(ne^{-t})}\dt}\\ 
		&&\Eq\  \Bl\hPhi(n)\sqln\Br^{-1}Z_{\log n}
		   \{\hPhi(n^{1-u(n)}) - \hPhi(n^{1-u})\} \\
	 &&\approx\ \Blb 1 - \frac{\hPhi(n^{1-u})}{\hPhi(n)}\Brb\,\s W_{\log n} (1),
\ens
whose randomness is determined only by the value of 
$W_{\log n}(1) \sim -(\t_n - \log n)/\s\sqln$.  To
match this with the corresponding formula for~$Y_n\uh(t)$, note
that, under~\Ref{L-cond2}, the second term in $G_n[U_n]$,
\[  
  \Bl \hPhi(n)\sqln \Br^{-1}
	  \int_{(\log n-t)_+}^{(\log n - t + \s U_n\sqln)_+} \hPhi(e^v)\dv\,,
\]
is small for $\log n - t = O(\sqln)$, and that, for larger values
of~$\log n - t = (1-u)\log n$, one can replace~$\hPhi(e^v)$	by
$\hPhi(n^{1-u})$ in the integral.  		 
\vskip0.5cm

\noindent
{\bf Remark.} Setting formally $\Phi(n)={\rm const}$ in the above formulas 
suggests that that $K_n\sim\tau_n$ in the case of bounded $\nu_0$. The latter is indeed true
and, moreover, 
$|K_n-\tau_n|$ remains bounded with all moments as $n$ grows;
the reason for this behaviour in the compound Poisson case is just that essentially
all gaps within ${\cal R}\cap [0,\log n]$ are hit by  the atoms of $Y_n$,
hence $K_n$ is close to the number of renewals on $[0,\log n]$.


\begin{thebibliography}{99}

\bibitem{ABT} R. Arratia, A.D. Barbour and S. Tavar{\'e}. {\it Logarithmic Combinatorial Structures:
A Probabilistic Approach}, European Math. Soc. Monographs in Math., v. 1, 2003.

\bibitem{Bertoin} J. Bertoin, {\it L{\'e}vy Processes}, Cambridge University Press, 1996.

\bibitem{BertoinSub} J. Bertoin, {\it Subordinators: Examples and Applications},
Springer Lecture Notes in Math. vol. 1727, 1996


\bibitem{BGT} N.H. Bingham, C.M. Goldie and J.L. Teugels. {\it Regular Variation}, 
Cambridge University Press,
1987.


\bibitem{Dudley} R.M. Dudley. 
{\it Probabilities and metrics\/},
Lecture Notes Series {\bf 45}, Aarhus Universitet, 1976.



\bibitem{Feller} W. Feller, {\it An Introduction to Probability Theory and its Applications}, volume II,
Wiley, 2nd edition, 1971.


\bibitem{Bernoulli} A.V. Gnedin. The Bernoulli sieve, {\it Bernoulli} 10: 79--96, 2004.

\bibitem{RCS} A.V. Gnedin and J. Pitman. Regenerative composition structures,
{\it Ann. Probab.} 33: 445--479, 2005.

\bibitem{RPS} A.V. Gnedin and J. Pitman. Regenerative partition structures,
{\it Elec. J. Combin.} 11(2),  paper \#R12, 2005.


\bibitem{GPYI} A.V. Gnedin, J. Pitman and M. Yor. Asymptotic laws for compositions derived from transformed
subordinators, {\it Ann. Probab.} (to appear), http:// arxiv.org/abs/math.PR/0403438, 2004.


\bibitem{GPYII} A.V. Gnedin, J. Pitman and M. Yor. Asymptotic laws for regenerative compositions:
gamma subordinators and the like, http:// arxiv.org/abs/math.PR/0405440, 2004.


\bibitem{LastBrandt} G. Last and A. Brandt. {\it Marked Point Processes on the Real Line: The Dynamic Approach},  
Springer, NY, 1995.



\bibitem{Mueller} D.W. M\"uller. Verteilungs-Invarianzprinzipien f\"ur das
starke Gesetz der gro\ss en Zahl, {\it Z.~Wahrscheinlichkeitstheorie verw.\ Geb.\/}
10: 173--192, 1968.

\bibitem{CSP} J. Pitman, {\it Combinatorial Stochastic Processes}, Springer Lecture Notes Math. (to appear),
 2002.


\bibitem{Winkel} M. Winkel, Electronic foreign exchange markets and level passage events of multivariate subordinators,
{\it J. Appl. Prob.} (to appear)


\end{thebibliography}
\end{document}